\documentclass[10pt]{amsart}
\usepackage[all]{xy}
\usepackage{graphicx}

\newcommand{\invliml}{\lower4pt\hbox{$\displaystyle{\varprojlim}\atop
\ell$}}
\newcommand{\dirlimk}{\lower4pt\hbox{$\displaystyle{\varinjlim}\atop
k$}}
\newcommand{\dirliml}{\lower4pt\hbox{$\displaystyle{\varinjlim}\atop\el
l$}}

\newtheorem{theo}{Theorem}[section]
\newtheorem{cor}[theo]{Corollary}
\newtheorem{prop}[theo]{Proposition}
\newtheorem{rem}[theo]{Remark}
\newtheorem{lem}[theo]{Lemma}
\newtheorem{defi}[theo]{Definition}
\newtheorem{ex}[theo]{Example}

\def\ZZ{\mathbb Z}
\def\RR{\mathbb R}

\begin{document}

\title{Rigidity and gluing for Morse and Novikov Complexes}

\author{Octav Cornea and Andrew Ranicki}

\address{O.C.: University of Montr\'eal\newline
\indent Department of Mathematics and Statistics
\newline \indent CP 6128 Succ. Centre Ville
\newline \indent Montr\'eal, QC H3C 3J7
\newline\indent Canada}
\email{cornea@dms.umontreal.ca}

\address{A.R.: Department of Mathematics and Statistics\newline
\indent University of Edinburgh\newline \indent
King's Buildings, Mayfield Road\newline \indent
Edinburgh EH9 3JZ\newline \indent Scotland, UK}
\email{aar@maths.ed.ac.uk}

\begin{abstract}\hskip5pt
We obtain rigidity and gluing results for the Morse complex of a
real-valued Morse function as well as for the Novikov complex of a
circle-valued Morse function.
A rigidity result is also proved for the Floer complex of a
hamiltonian defined on a closed symplectic manifold $(M,\omega)$
with $c_{1}|_{\pi_{2}(M)}=[\omega]|_{\pi_{2}(M)}=0$.
The rigidity results for these complexes
show that the complex of a fixed generic
function/hamiltonian is a retract of the Morse (respectively Novikov or
Floer) complex of
any other sufficiently $C^{0}$ close generic function/hamiltonian.
The gluing result is a type of Mayer-Vietoris formula for the Morse
complex.
It is used to express algebraically the Novikov complex up to
isomorphism in terms
of the Morse complex of a fundamental domain. Morse cobordisms are used
to compare various Morse-type complexes without the need of bifurcation
theory.
\end{abstract}
\maketitle

\tableofcontents

\section*{Introduction.}
A Morse type complex associated to a generic pair $(f,\alpha)$
with $f$ a functional and with $\alpha$ an additional geometric
structure is a complex freely generated by the critical points
(assumed in the paper to be finite in number) of $f$ and with a
differential given by counting ``flow lines" that join successive
critical points of $f$. The definition of these ``flow lines''
depends on both $f$ and $\alpha$. For example, $f$ might be a
Morse function, $\alpha$ a Riemannian metric and then the
differential counts flow lines of $-\nabla^{\alpha}f$. The complex
is in this case the classical Morse complex variously attributed
to Thom, Smale, Milnor and under a different form to Witten. Other
examples that will appear in this paper are provided by the
Novikov complex where $f$ is an $S^{1}$-valued function (the
additional structure being again a metric) and the Floer complex
which, from this perspective, is associated to an action
functional defined on a space of unbased, contractible loops on a
symplectic manifold together with an almost complex structure on
the manifold. In this last example, the ``flow lines" to be
counted are  deformed pseudo-holomorphic cylinders. In this
introduction we shall denote any complex of one of these three
types by $\mathcal{C}(f,\alpha)$.

\

Our focus in this paper is not on the homology computed by
$\mathcal{C}(f,\alpha)$ (which, indeed,
is in all these cases independent of the pair $(f,\alpha)$) but rather
on the complexes themselves.

\

The first type of result that we prove for the Morse, Novikov as well
as for the Floer complex
is a rigidity statement:  it shows that small $C^{0}$ perturbations of
$f$
accompanied by arbitrary perturbations of $\alpha$ do not decrease the
complexity  of
$\mathcal{C}(f,\alpha)$.
Precisely, for fixed $(f,\alpha)$ there exists
$\delta_{f}>0$ such that if $(f',\alpha')$ is another generic pair such
that
$\mathcal{C}(f',\alpha')$ is defined and with $||f-f'||_{0}\leq
\delta_{f}$, then $\mathcal{C}(f,\alpha)$
is a retract of $\mathcal{C}(f',\alpha')$ (which means that there are
chain morphisms
$i:\mathcal{C}(f,\alpha)\to\mathcal{C}(f',\alpha')$,
$j:\mathcal{C}(f',\alpha')\to\mathcal{C}(f,\alpha)$
such that $j\circ i=id$). When $f=f'$ we immediately deduce
that the isomorphism type of $\mathcal{C}(f,\alpha)$ is independent of
$\alpha$.

\

For the classical Morse complex and for the Novikov complex,
this corollary is implicit  in the work of Latour \cite{Lat}
where the approach is that of bifurcation analysis.
The general case, when $f'$ is different from (but $C^{0}$-close to)
$f$,
is new even in the standard Morse case and bifurcation analysis does
not suffice to
prove it. The reason is that $C^{0}$-closeness to $f$ is not naturally
encoded in
some restriction of the type of modification that may occur in the
Morse complex at
passage through bifurcation points when tracing a generic path of
(of $C^{\infty}$ functions) from $f$ to $f'$. Of course, without the
$C^{0}$-closeness condition the result fails. Moreover, in the Floer
setting bifurcation
theory is very little understood and therefore particularly inefficient
for this type of
problem.

\

In this paper we do not use bifurcation analysis at all. Rather,
we introduce and apply a notion of Morse cobordism. Roughly
speaking, this is a generic pair $(f,\alpha)$ with $f$ a Morse
function and $\alpha$ a Riemannian metric  both defined on a
cobordism $(M;N_{0},N_{1})$ and such that $\nabla^{\alpha}f$ is
tangent to the boundary and the restrictions $f_{i}=f|_{N_{i}}$
are Morse functions whose critical points verify
$ind_{f_{i}}(x)=ind_{f}(x)+i$, $i\in\{0,1\}$. The role of Morse
cobordisms is to transform a $k$-parametric deformation question
for Morse functions or complexes in a problem concerning a Morse
function (or complex) defined on a space with $k$ more dimensions
(in the paper the actual values of $k$ that are used are just $1$
and $2$). We introduce this notion in Section \ref{Morse
cobordism} and discuss its immediate properties in \S \ref{basic}.
In \S \ref{subsubsec:rigid} we show the rigidity result for the
Morse complex. This serves as prototype for the  Floer and Novikov
cases which are shown respectively in \S\ref{Floer_rig} and \S
\ref{Novikov} Theorem \ref{theo:Novikov} (i). It should be noted
that the whole idea of Morse cobordisms is inspired by the
continuation method in Floer theory. Indeed, the special Morse
cobordisms without critical points in the interior of $M$ (which
implies that $(M;N_{0},N_{1})$ is a trivial cobordism) have an
excellent Floer theory translation provided by the {\em monotone}
homotopies introduced by Floer and Hofer in \cite{FlHo}  (see also
\cite{CiFlHo}) and which we shall use in \S\ref{Floer_rig}.

\

In \S \ref{gluing} we show a Mayer-Vietoris type formula for the Morse
complex.
A geometric cobordism $(M;N_{0},N_{1})$ can be split along a two-sided
hypersurface $N\subset Int(M)$
into two cobordisms $(M';N_{0},N)$, $(M'';N,N_{1})$. Of course, the
initial cobordism can easily be
recovered by gluing back the two pieces. We discuss how to split a
Morse cobordism
$f:(M;N_{0},N_{1})\to \RR$ along a Morse function $g:N\to \RR$ into two
Morse
cobordisms $h':M'\to\RR$ and $h'':M''\to\RR$. The
key point showed in the gluing Theorem \ref{theo:glueadapt},
is that, if a certain technical property is satisfied by $g$, then it
is possible to
express in a purely algebraic manner the Morse complex of $f$ in terms
of the Morse complexes of $h'$ and
$h''$. Again, for this result bifurcation methods are not directly
usable
(see Remark \ref{rem:glue1}).
Even if $g$ does not satisfy the technical property indicated
the same method still produces the isomorphism type of the Morse
complex of $f$.

\

The result in \S\ref{gluing} is applied to the study of the Novikov
complex in \ref{Novikov}. In Theorem \ref{theo:Novikov} (ii-iii) we
describe
a way to recover this complex algebraically out of information on one
fundamental domain.
For a fixed Morse-Smale function $(f,\alpha):M\to S^{1}$ the procedure
consists
in first splitting $f$ along a Morse function $g:N\to \RR$ with $N$ a
regular level
surface of $f$. This produces a Morse cobordism on one fundamental
domain of $f$.
We then consider the Morse complex associated to this Morse cobordism
and,
by a purely algebraic procedure, we glue together an infinite
number of copies of it thus getting a new complex $\mathcal{C}$. The
problem then
is to compare $\mathcal{C}$ with $\mathcal{C}(f,\alpha)$. By using
again
Morse cobordisms we show that these two complexes are isomorphic.
The complex $\mathcal{C}$ is viewed as an approximation of
$\mathcal{C}(f,\alpha)$
and special  choices of $g$ result in a complex $\mathcal{C}$
agreeing with the Novikov complex of $f$ with any finite
degree of precision desired. This last point is related to
results of Pajitnov \cite{P} which provide an exact formula for the
Novikov
complex but only apply to a special class of Morse-Smale functions.
We discuss these relations in Remark \ref{rem:fin1}(c).

\section{Cobordism of Morse functions}\label{Morse cobordism}

\subsection{Basic constructions}\label{basic}

This subsection is concerned with the basic definitions and properties
of Morse cobordisms.

\subsubsection{Definitions.} We start by introducing an extension of
the usual notion of a
Morse-Smale function on a closed manifold to a manifold with boundary.
Our definition is designed to ensure that the standard construction of
the chain complex in the closed case extends to a construction of a
chain complex for a Morse-Smale function on a manifold with boundary.

Let $M$ be a compact manifold, possibly with non-empty boundary.  A
{\it Morse function} $f:M \to \RR$ is a smooth function with
nondegenerate critical points (that might belong to the boundary).  We
denote by ${\rm Crit}_i(f)$ (or ${\rm Crit}_i(M,f)$) the critical
points $x$ of index $i$, ${\rm ind}_f(x)=i$, of such a function and we
let ${\rm Crit}_{\ast}(f)$ (or ${\rm Crit}_{\ast}(M,f)$) be the set of
all critical points of $f$.

Given a Morse function $f:M \to \RR$ and a Riemannian metric $\alpha$
on $M$ let $\nabla f$ be the gradient vector field of $f$.  Let
$\phi=\{\phi_t:M \to M : t \in \RR\}$ be the flow (possibly only
partially defined) on $M$ induced by $-\nabla f$, so that
$$d \phi_t(x)/d t ~=~-\nabla f (\phi_t(x))~~,~~
\phi_0(x)~=~x~~(x \in M)~.$$
For a critical point $p\in M$ of $f$ let
$$W_f^u(p)~=~\{x\in M : \lim_{t\rightarrow -\infty}\phi_t(x)=p\}$$
be the unstable manifold of $p$ and let
$$W_f^s(p) ~=~\{x\in M : \lim_{t\rightarrow\infty}\phi_t(x)=p\}$$
be the stable manifold of $p$.

\begin{defi}\label{Morse-Smale} {\rm A {\it Morse-Smale function}
$(f,\alpha)$ on a compact manifold with boundary $(M,\partial M)$ is a
Morse function $f:M\to \RR$ together with a Riemannian metric $\alpha$
such that
\begin{itemize}
\item[(i)] for each boundary component $N$ of
$\partial M$ one of the two following conditions is satisfied:
\begin{itemize}
\item[(a)] $f$ is regular on $N$, meaning that for all
$x \in N$ $$\nabla f(x) \notin T_xN \subset T_xM~.$$
\item[(b)] For all $x\in N$ we have $$\nabla f(x) \in T_xN~.$$
\end{itemize}
\item[(ii)] for any two critical points $p,q\in M$ the stable and
unstable manifolds $W^u(p)$ and $W^s(q)$ intersect transversely (due to
point (i) the stable and unstable manifolds are indeed submanifolds of
$M$, possibly with boundary if $\partial M\not=\emptyset$).
\hfill$\qed$
\end{itemize}}
\end{defi}

\begin{rem} {\rm The classical definition of a Morse-Smale function
assumes the boundary behaviour of (i) (a).  Let
$(g,\beta)=(f,\alpha)|_N$ be defined when condition (i) (b) is
satisfied.  Then $(g,\beta)$ is also a Morse-Smale function on $N$.  In
this case, ${\rm Crit}_{\ast}(g)\subset {\rm Crit}_{\ast}(f)$ (because
the component of $\nabla f(x)$ normal to the boundary is null for all
$x\in N$) and if $x\in {\rm Crit}_{\ast}(g)$ we have
$${\rm ind}_g(x)~=~{\rm ind}_f(x)-\delta_x$$
with $\delta_x\in\{0,1\}$.\hfill\qed}
\end{rem}

Let $(f,\alpha):M \to \RR$ be Morse-Smale.  Denote the set of critical
points of $f$ with index $i$ by ${\rm Crit}_i(f)$.  For $p\in {\rm
Crit}_i(f)$, $q\in {\rm Crit}_j(f)$ the space $Z(p,q)$ of (non-broken)
flow lines of $\phi$ that join $p$ to $q$ is homeomorphic to any
intersection $W^u(p)\cap W^s(q)\cap f^{-1}(a)$ where $a$ is some
regular value of $f$ such that $f(q)<a<f(p)$.  The space
$$Z(p,q)~\cong~W^u(p)\cap W^s(q)\cap f^{-1}(a)$$
is a $(i-j-1)$-dimensional manifold.  If $j=i-1$ (and as $M$ is
compact) this is a finite set.  We fix orientations of $T_xM$ for all
$x\in {\rm Crit}_{\ast}(f)$ as well as for each tangent space to the
unstable manifolds $T_x(W^u(x))$.  This provides orientations for
$T_x(W^s(x))$ by demanding that the orientation on $T_x(W^u(x))$ and
that on $T_x(W^s(x))$ give that fixed on $T_x(M)$.  Moreover, as
$W^s(x)$, $W^u(x)$ are contractible these orientations induce
orientations of the tangent spaces to the whole stable and unstable
manifolds.  For a point $z\in Z(p,q)$ we now let $\epsilon(z)$ be equal
to $1$ if the orientation of $T_z(W^u(p))$ and that of $T_z(W^s(q))$
(in this order) give the orientation induced on $T_zM$ from $T_p(M)$.
Finally, let
$$n^{f,\alpha}(p,q)~=~\sum\limits_{z\in Z(p,q)}\epsilon(z) \in \ZZ~.$$

\begin{defi}\label{Morse complex1}
{\rm The {\it $\ZZ$-coefficient Morse complex} $C(M,f,\alpha)$ of a
Morse-Smale function $(f,\alpha):M \to \RR$ is defined by
$$\begin{array}{l}
d~:~C(M,f,\alpha)_i~=~\ZZ[{\rm Crit}_i(f)]
\to C(M,f,\alpha)_{i-1}~=~\ZZ[{\rm Crit}_{i-1}(f)]~;\\[1ex]
\hskip100pt p \mapsto \sum\limits_{q\in {\rm
Crit}_{i-1}(f)}n^{f,\alpha}(p,q)q~.
\end{array}$$
\hfill\qed}
\end{defi}

The $\ZZ[\pi]$-coefficient Morse complex is defined for a Morse-Smale
function $(f,\alpha):M \to \RR$ and a regular cover $\widetilde{M}$ of
$M$ with group of covering translations $\pi$.  Let
$(\widetilde{f},\widetilde{\alpha}):\widetilde{M} \to \RR$ be the
pullback of $(f,\alpha)$, and let $\widetilde{\phi}$ be the pullback of
$\phi$.  The critical points $\widetilde{p} \in \widetilde{M}$ of
$\widetilde{f}$ are the lifts of the critical points $p \in M$ of $f$.
Fix two critical points $\widetilde{p} \in {\rm Crit}_i(\widetilde{f})$
and $\widetilde{q} \in {\rm Crit}_j(\widetilde{f})$.  As each path in
$M$ that originates at $p$ lifts to a unique path in $\widetilde{M}$ of
origin $\widetilde{p}$, the space $\bigcup_{g\in\pi}
Z(\widetilde{p},g\widetilde{q})$ of flow lines of $\widetilde{\phi}$
that join $\widetilde{p}$ to one of the points $g\widetilde{q}$,
$g\in\pi$ is homeomorphic to $Z(p,q)$.  In particular, for $j=i-1$ the
sum below is well defined and satisfies :
$$\sum\limits_{g\in \pi}
n^{\widetilde{f},\widetilde{\alpha}}(\widetilde{p},
g\widetilde{q})~=~ n^{f,\alpha}(p,q) \in \ZZ$$
and
$$n^{\widetilde{f},\widetilde{\alpha}}(g\widetilde{p},g\widetilde{q})~=
~
n^{\widetilde{f},\widetilde{\alpha}}(\widetilde{p},\widetilde{q})
\in \ZZ~~ (g \in \pi)~.$$
Choose a lift of each critical point $p \in {\rm Crit}_i(f)$ to a
critical point $\widetilde{p} \in {\rm Crit}_i(\widetilde{f})$,
allowing the identification
$${\rm Crit}_i(\widetilde{f})~=~\pi \times {\rm Crit}_i(f)~.$$

\begin{defi}\label{Morse complex2}
{\rm The {\it $\ZZ[\pi]$-coefficient Morse complex} $C(M,f,\alpha)$ is
the
based f.g. free $\ZZ[\pi]$-module chain complex given by
$$\begin{array}{l}
d~:~C(M,f,\alpha)_i~=~\ZZ[\pi][{\rm
Crit}_i(f)]\to
C(M,f,\alpha)_{i-1}~=~\ZZ[\pi][{\rm Crit}_{i-1}(f)]~;\\[1ex]
\hskip150pt
p \mapsto \sum\limits_{q\in {\rm Crit}_{i-1}(f)}(\sum\limits_{g \in
\pi}
n^{\widetilde{f},\widetilde{\alpha}}(\widetilde{p},g\widetilde{q})g)q~.
\end{array}$$
\hfill\qed}
\end{defi}

This clearly depends on the choices of the lifts $\widetilde{p}$. More
 invariantly, the $\ZZ[\pi]$-coefficient Morse complex can be written
as
$$\begin{array}{l}
d~:~C(M,f,\alpha)_i~=~\ZZ[{\rm
Crit}_i(\widetilde{f})]\to
C(M,f,\alpha)_{i-1}~=~\ZZ[{\rm Crit}_{i-1}(\widetilde{f})]~;\\[1ex]
\hskip150pt \widetilde{s} \mapsto
\sum\limits_{\widetilde{q}\in {\rm Crit}_{i-1}(\widetilde{f})}
n^{\widetilde{f},\widetilde{\alpha}}(\widetilde{s},\widetilde{q})
\widetilde{q}~.
\end{array}$$

\begin{rem}\label{rem:class} {\rm (a) In the classical situation
$$(f,\alpha):(M;N_0,N_1) \to ([0,1];\{0\},\{1\})$$
is a Morse-Smale function on a cobordism (Milnor \cite{M}).  This fits
into our definition of a Morse-Smale function (\ref{Morse-Smale})
provided that $f$ is regular on $N_i$ and its negative gradient points
"out" on $N_0$ and "in" on $N_1$.  Such a Morse-Smale function on an
$m$-dimensional cobordism $(M;N_0,N_1)$ determines a handle
decomposition
$$M~=~N_0 \times [0,1] \cup \bigcup^m_{i=0}\bigcup_{{\rm Crit}_i(f)}D^i
\times D^{m-i}~.$$
The handle decomposition of $(M;N_0,N_1)$ gives $(M,N_0)$ the structure
of a relative $CW$ pair, with one $i$-cell for each $i$-handle.
Franks' paper \cite{F} identifies the Morse complex of $f$ with the
associated cellular chain complex $C(M;\ZZ[\pi])=C(\widetilde{M})$.  It
is easy to verify that the definition gives a chain complex also in our
more general context.\\
(b) Morse-Smale functions $f:M\to \RR$ have the important property that
if $p\in {\rm Crit}_{\ast}(f)$ and $\{x_n\}\subset W^{u}(p)$ is a
convergent sequence in $M$ with limit $x_{\infty}\in M$, then
$x_{\infty}$ is situated on a possibly broken flow line originating at
$p$ (that is a flow line that passes geometrically through some other
critical points besides $p$ before arriving in $x_{\infty}$).  In fact,
one of the proofs that the Morse complex is indeed a complex is
obtained by understanding precisely the natural compactifications of
the spaces of flow lines.  \hfill\qed}
\end{rem}

To simplify notation, we shall write $C(M,f)$ (resp.  $n^f(p,q)$)
instead of $C(M,f,\alpha)$ (resp.  $n^{f,\alpha}(p,q)$) whenever the
choice of metric $\alpha$ is clear from the context.  Also, by an abuse
of terminology, we shall write $C(M;\ZZ[\pi])=C(\widetilde{M})$ as
$C(M)$, and $H_*(M;\ZZ[\pi])=H_*(\widetilde{M})$ as $H_*(M)$.

\begin{defi}\label{cobordism}
{\rm A {\it cobordism} from a Morse-Smale function
$$(g_1,\beta_1)~:~N_1 \to \RR$$ to a Morse-Smale function
$$(g_0,\beta_0)~:~N_0 \to \RR$$
is a cobordism $(M;N_0,N_1)$ together with a
Morse-Smale function
$$(f,\alpha)~:~M \to \RR$$
such that
\begin{itemize}
\item[(i)] $(f,\alpha)|_{N_0}= (g_0,\beta_0)$ with
$${\rm ind}_f(p)~=~{\rm ind}_{g_0}(p)$$
for each critical point $p\in N_0$,  \\[1ex]
\item[(ii)] $(f,\alpha)|_{N_1}= (g_1,\beta_1)$ with
$${\rm ind}_f(p)~=~{\rm ind}_{g_1}(p)+1$$
for each critical point $p\in N_1$.
\end{itemize}}
\hfill\qed
\end{defi}
\vskip-5mm
$$\xymatrix@C+20pt{\ar@{-}[rrrr] \ar@{-}[dd]_{\displaystyle{N_0}} &&
\ar@{}[dd]^{\hskip-2.5mm\displaystyle{M}} && \ar@{-
}[dd]^{\displaystyle{N_1}}\\
&&&&\\
\ar@{-}[rrrr]
\ar[dd]_{\displaystyle{(g_0,\beta_0)}} &&
\ar[dd]^{\displaystyle{(f,\alpha)}} &&
\ar[dd]^{\displaystyle{(g_1,\beta_1)}}\\
&&&&\\
\ar@{-}[rrrr]_{\displaystyle{\RR}} &&&&}
$$
\medskip

\begin{rem} {\rm
(a) We denote a Morse cobordism as before by
$$(f,\alpha):(M;N_0,N_1)\to \RR~.$$
Clearly, $N_0$ and $N_1$ play different roles in the definition above.
However,
because we only work here with compact manifolds, it follows from Lemma
\ref{constr2} below that Morse-Smale functions $(g_0,\beta_0):N_0 \to
\RR$,
$(g_1,\beta_1):N_1 \to \RR$ are cobordant if and only if the manifolds
$N_0$, $N_1$
are cobordant and therefore the Morse-cobordism relation is an
equivalence.\\
(b) The condition (i) in Definition \ref{cobordism} together with the
Definition
\ref{Morse-Smale} imply that no negative gradient flow lines of
$(f,\alpha)$ can
leave $N_{0}$ towards $M\backslash \partial M$. Similarly, no
trajectory can enter
$N_{1}$ from the interior due to condition (ii). } \hfill\qed
\end{rem}

\begin{ex} {\rm
An $m$-dimensional cobordism $(M;N_0,N_1)$ admits an embedding
$$(M;N_0,N_1) \subset \RR^n$$
for $n \geq 2m+1$, with
$$\begin{array}{l}
M \subset \{(x_1,x_2,\dots,x_n) \in \RR^n\,:\,0 \leq x_n \leq
1\}~,\\[1ex]
N_0~=~\{(x_1,x_2,\dots,x_n) \in M\,:\,x_n=0\}~,\\[1ex]
N_1~=~\{(x_1,x_2,\dots,x_n) \in M\,:\,x_n=1\}
\end{array}$$
and such that the height function
$$f~:~M \to \RR~;~(x_1,x_2,\dots,x_n) \mapsto x_n$$
is Morse-Smale on $M\backslash \partial M$.  As $f$ is constant on
$N_0$ and $N_1$, it is not a cobordism of Morse-Smale functions.
However, it is possible to slightly tilt $M$ inside $\RR^n$ to obtain
an isotopic embedding $(M;N_0,N_1) \subset \RR^n$ such that the height
function $M \to \RR$ restricts to Morse-Smale functions on $N_0$ and
$N_1$.  This is not yet a cobordism of Morse-Smale functions as in
\ref{cobordism} because the component of the gradient of the height
function that is normal to $\partial M$ does not generally vanish.  As
we shall see in Lemma \ref{constr2} below, it is possible to perturb
the height function in a neighbourhood of $\partial M$ such that the
resulting function becomes a cobordism of Morse-Smale functions.}
\hfill\qed \end{ex}

\subsubsection{Immediate Properties of Morse cobordisms.}

We recall that the {\it algebraic mapping cone}
$C(\phi)$ of a chain map $\phi:C
\to D$ is the chain complex defined by
$$d_{C(\phi)}~=~\begin{pmatrix}
d_D & \phi \\
0 & -d_C \end{pmatrix}~:~ C(\phi)_i~=~D_i \oplus C_{i-1} \to
C(\phi)_{i-1}~=~D_{i-1} \oplus C_{i-2}~.$$
A chain homotopy $\psi:\phi\simeq\phi':C\to D$ determines an
isomorphism of the algebraic mapping cones
$$h~=~\begin{pmatrix} 1 & \psi \\
0 & 1 \end{pmatrix}~:~C(\phi)\to C(\phi')~.$$
If $C,D$ are based f.g.  free $\ZZ[\pi]$-module chain complexes (as
will be the case in the applications to Morse complexes) then $h$ is a
simple isomorphism of based f.g.  free $\ZZ[\pi]$-module chain
complexes.

The Morse complex of a cobordism of Morse-Smale functions has the
following homological properties :

\begin{prop} \label{chain complex}
Let $(f,\alpha):(M;N,N') \to \RR$ be a
cobordism of Morse-Smale functions $(g,\beta):N
\to \RR$, $(g',\beta'):N' \to \RR$, write
$$D~=~C(N,g)~~,~~D'~=~C(N',g')~~,~~F~=~C(M\backslash\partial
M,f\vert)~,$$
and let $\widetilde{M}$ be a regular
cover of $M$ with group of covering translations $\pi$. \\
{\rm (i)} There are 3 types of critical points
of $f$
$${\rm Crit}_i(M,f)~=~{\rm Crit}_i(N,g)\cup
{\rm Crit}_i(M\backslash \partial M,f\vert)\cup
{\rm Crit}_{i-1}(N',g')~.$$ The Morse complex
of $f$ is given by
$$\begin{array}{l}
C(M,f)_i~=~\ZZ[\pi][{\rm Crit}_i(M,f)]~=~D_i \oplus F_i\oplus D'_{i-
1}~,\\[1ex]
d_{C(M,f)}~=~\begin{pmatrix}
d_D & \theta & -\psi \\
 0  & d_F & -\theta' \\
 0  & 0 & -d_{D'} \end{pmatrix}~:\\[2ex]
C(M,f)_i~=~D_i \oplus F_i \oplus D'_{i-1} \to
C(M,f)_{i-1}~=~D_{i-1} \oplus F_{i-1} \oplus
D'_{i-2}
\end{array}$$
with
$$\begin{array}{l}
\theta~:~F_i \to D_{i-1}~;~
p \mapsto \sum\limits_{q \in {\rm Crit}_{i-
1}(g)}n^{f,\alpha}(p,q)q~,\\[1ex]
\theta'~:~D'_i \to F_i~;~
p \mapsto \sum\limits_{q \in {\rm
Crit}_i(e)}n^{f,\alpha}(p,q)q~,\\[1ex]
\psi~:~D'_i \to D_i~;~ p \mapsto \sum\limits_{q
\in {\rm Crit}_i(g)}n^{f,\alpha}(p,q)q
\end{array}$$
such that
$$\begin{array}{l}
d_D\theta + \theta d_F~=~0~:~F_i \to D_{i-2}~,\\[1ex]
d_F\theta' - \theta' d_{D'}~=~0~:~D'_i \to F_{i-1}~,\\[1ex]
d_D\psi+\psi d_{D'}+\theta\theta'~=~0~:~D'_i
\to D_{i-1}~.
\end{array}$$
{\rm (ii)} The Morse complex $C(M,f)$ is simple chain equivalent to
$C(M,N')$, the $\ZZ[\pi]$-coefficient cellular complex for any relative
$CW$ structure on $(M,N')$, with homology
$$H_*(C(M,f))~=~H_*(M,N')~.$$
{\rm (iii)} The subcomplex of $C(M,f)$
$$C(M\backslash N',f\vert)~=~C(\theta)~=~\big(~D_i \oplus F_i~,~
\begin{pmatrix} d_D & \theta \\  0  & d_F \end{pmatrix}~\big)$$
is simple chain equivalent to $C(M)$, with homology
$$H_*(C(M\backslash N',f\vert))~=~H_*(M)~.$$
{\rm (iv)} The quotient complex of $C(M,f)$
$$C(M\backslash N,f\vert)~=~C(\theta')~=~
\big(~F_i \oplus D'_{i-1}~,~\begin{pmatrix}
d_F & -\theta'  \\
0 & -d_{D'} \end{pmatrix}~\big)$$ is simple
chain equivalent to $C(M,\partial M)$, with
homology
$$H_*(C(M\backslash N,f\vert))~=~H_*(M,\partial M)~.$$
{\rm (v)} The subquotient complex of $C(M,f)$
$$C(M\backslash \partial M,f\vert)~=~F$$
is simple chain equivalent to $C(M,N)$, with
homology
$$H_*(C(M\backslash \partial M,f\vert))~=~H_*(F)~=~H_*(M,N)~.$$
{\rm (vi)} The chain maps (up to sign)
$$\begin{array}{l}
\theta~:~F~=~C(M\backslash \partial M,f\vert) \to D_{*-1}~=~C(N,g)_{*-
1}~,\\[1ex]
\theta'~:~D'~=~C(N',g') \to F~=~C(M\backslash
\partial M,f\vert)
\end{array}$$
induce the natural morphisms in homology
$$\begin{array}{l}
\theta_*~=~\partial~:~H_*(F)~=~H_*(M,N) \to H_{*-1}(D)~=~H_{*-
1}(N)~,\\[1ex]
\theta'_*~:~H_*(D')~=~H_*(N') \to
H_*(F)~=~H_*(M,N)
\end{array}$$
and $\psi$ is a chain homotopy
$$\psi~:~\theta\theta'~\simeq~0~:~D'~=~C(N',g') \to D_{*-
1}~=~C(N,g)_{*-1}~.$$
\end{prop}

\begin{proof} (i) For any $x \in N$
$$\lim_{t\rightarrow\infty}\phi_t(x) \in N~,$$
corresponding to the entry $d_D$ in the first column of $d_{C(M,f)}$.\\
For any $x \in M\backslash \partial M$
$$\lim_{t\rightarrow\infty}\phi_t(x) \in N~{\rm or}~M\backslash
\partial M~,$$
corresponding to the two entries $\theta,d_F$
in the second column of
$d_{C(M,f)}$.\\
For any $x \in N'$
$$\lim_{t\rightarrow\infty}\phi_t(x) \in N'~{\rm or}~M\backslash
\partial M~
{\rm or}~N~,$$ corresponding to the three
entries $-\psi,-\theta',-d_{D'}$
in the third column of $d_{C(M,f)}$.\\
(ii) Slightly extend $M$ to a manifold $M'$ by pasting to $M$ collars
homeomorphic to $N\times [0,1]$ and $N' \times [0,1]$.  We also extend
the function $f$ and the metric $\alpha$ to a function $f'$
respectively a metric $\alpha'$ in such a way that $f'$ has the same
critical points as $f$ and is regular and its gradient points out on
$N'\times
\{1\}$ and points inside and is regular on $N\times \{1\}$.  Then, by
standard Morse theory, $C(M,f',\alpha')$ is simple chain equivalent to
$C(M',N')$.  On the other hand we have
$C(M,f',\alpha')=C(M,f,\alpha)$.\\
(iii) The projection from $C(M,f)$ to $C(N',g')_{*-1}$ is given by
$$C(M,f)~\simeq~C(M,N')  \xymatrix{\ar[r]^-{\displaystyle{\partial}}&}
C(N',g')_{*-1}~\simeq~C(N')_{*-1}~.$$ The kernel of this projection is
precisely
$C(M\backslash N', f|)$ and this complex is therefore simple chain
equivalent to
$C(M)$.\\
(iv) The chain inclusion
$$C(N,g)~\simeq~C(N)
\to C(M,f)~\simeq~C(M,N')$$
is a chain representative of the inclusion $N\hookrightarrow (M,N')$.
Therefore, $C(M\backslash N, f|)$ (which is the co-kernel of the map
above) is simple chain equivalent to $C(M,N\cup N')=C(M,\partial M)$.\\
(v) The chain inclusion
$$C(N,g)~\simeq~C(N) \to C(M\backslash N',f\vert)~\simeq~C(M)$$
is a chain representative of the inclusion $N\hookrightarrow M$.\\
(vi) Combine (ii), (iii), (iv) and (v).
\end{proof}


A particular type of Morse cobordism will play an important role
further on.

\begin{defi}\label{simple}
{\rm A cobordism of Morse-Smale functions
$$(f,\alpha)~:~(M;N,N')\to \RR$$
is {\it simple} if $f$ does not have any
critical points in $M\backslash \partial M$, so
that
$${\rm Crit}_i(M,f)~=~{\rm Crit}_i(N,g) \cup {\rm Crit}_{i-1}(N',g')~.
\eqno{\qed}$$}\end{defi}
\medskip

\begin{prop} \label{chaineq}
Let $(f,\alpha):(M;N,N') \to \RR$ be a simple cobordism of Morse-Smale
functions $(g,\beta):N \to \RR$, $(g',\beta'):N' \to \RR$.\\
{\rm (i)} The identity $N \to N$ extends to a diffeomorphism
$$(M;N,N') \to N \times ([0,1];\{0\},\{1\})~.$$
{\rm (ii)} The inclusions $N \hookrightarrow M$, $N' \hookrightarrow M$
induce simple chain equivalences
$$C(N,g)~\simeq~C(M)~~,~~C(N',g')~\simeq~C(M)$$
with $C(M)$ the $\ZZ[\pi]$-coefficient cellular
chain complex of any $CW$ structure on $M$.\newline
{\rm (iii)} The chain map
$$f^M~:~C(N',g')~\to~C(N,g)$$
defined by
$$\begin{array}{ll}
f^M~:~C(N',g')_i~=~\ZZ[\pi][{\rm
Crit}_i(g')]&\to
C(N,g)_i~=~\ZZ[\pi][{\rm Crit}_i(g)] ~;\\[1ex]
&x \mapsto \sum\limits_y n^f(x,y)y
\end{array}$$
is a simple chain equivalence, with $n^f(x,y)$
the algebraic number of downward gradient flow
lines in $\widetilde{M}$ from $x$ to $y$, for
any critical points $x \in \widetilde{N}'$, $y\in \widetilde{N}$ with
$${\rm ind}_{g'}(x)~=~{\rm ind}_f(x)-1~=~{\rm ind}_{g}(y)~=~
{\rm ind}_f(y)~=~i~.$$
\end{prop}
\begin{proof} (i) It is possible to find a small perturbation of the
negative gradient flow of $f$ in a neighbourhood of $\partial M$ such
that the
resulting flow $\gamma'$ points inside $M$ on $N'$ and outside on $N$;
it is
gradient like and has no stationary points inside $M$. The existence of
such a flow
implies the claim.  To construct this perturbation fix some Morse
charts $U_i$
around the critical points of $f$ that are situated on $N$.  Fix also a
collared
neighbourhood of $N$ that is diffeomorphic to $N\times [0,\epsilon]$.
With respect
to this parametrization notice that, for all small enough $\delta$ the
projection of
$\nabla f(x)$ on $T_{x}(N\times\{\delta\})$ is non-zero for all
sufficiently small
$\delta$ and $x\not \in \cup U_i$. Moreover, for $\delta$ small and
$x\in U_i$ we
have that $\nabla f(x)$ is not tangent to $T_{x}(N\times \{\delta\})$
and points
towards $N\times \{0\}=N$ (because of the Morse lemma).  This means
that by adding
to $-\nabla f$ a vector field which is $0$ outside $N\times [0,\tau]$
($\tau$ small
enough) and is equal to $-h(t)\partial/\partial t$ for $x=(z,t)\in
N\times [0,\tau]$
with $h:[0,\tau]\to [0,1]$ decreasing, $h(0)=1$, $h(\tau)=0$ the
induced flow will
satisfy the desired properties with respect to $N$.  Of course, a
similar
construction is possible relative to $N'$ and produces $\gamma'$.\\
(ii) Immediate from (i).\\
(iii) The chain map $f^M$ is the chain homotopy
of Proposition \ref{chain complex} (vi)
$$\psi~=~f^M~:~0~\simeq~0~:~ C(N',g') \to C(N,g)_{*-1}~.$$
For any $CW$ structure on $M$ the chain map
$f^M$ is chain homotopic to the composite of
the simple chain equivalences given by (ii)
$$C(N',g') \to C(M) \to C(N,g)~.$$
\end{proof}

\subsubsection{Construction of Morse-Smale
cobordisms.}\label{subsubsec:constr}

We shall discuss here a few simple ways to construct Morse cobordisms.
We first fix some more notation.
For $i=0,1,2,\dots$ let $C^i(N,\RR)$ be the space of $C^i$-functions
$h:N \to \RR$, with the following topology.  For $i=0$ use the norm
$$||h||_0~=~{\rm sup}\,\{\vert h(x)\vert : x\in N\}~.$$
For $i \geq 1$ use the Whitney $C^i$-topology, in which a neighbourhood
of $f \in C^i(N,\RR)$ consists of those $g\in C^i(N,\RR)$ such that in
local coordinates, $f$ and $g$ together with their first $i$
derivatives are within $\epsilon$ at each point of $N$ (Hirsch
\cite{Hi}, page 35).

Here is a result that summarizes the output of some of our
constructions.
All the arguments here are simple once the correct statements are
formulated -- similar constructions have been independently used by
Abbondandolo and Majer \cite{AM}.  We include the details because
we shall use these constructions repeatedly later in the paper.

\begin{prop}\label{cob_constr}
Fix a Morse-Smale function $(h,\gamma)$ on a compact closed manifold
$N$. Consider
$(h',\gamma'):N \to \RR$, another Morse-Smale function on $N$, such
that
$h(x)>h'(x)$, $\forall x\in N$. There exists a simple Morse-Smale
cobordism
$(H,\Gamma):N\times [0,1]\to \RR$ from $(h,\gamma)$  to $(h',\gamma')$
such that
$$h'(x)\leq H_t(x)\leq h(x)~~(x \in N,t\in [0,1])~.$$
The resulting chain map
$$i~=~H^{N\times [0,1]}~:~C(N,h,\gamma)\to  C(N,h',\gamma')$$
is a simple chain equivalence.
\end{prop}

\begin{proof} The proof is an immediate consequence of the
constructions of Morse cobordisms described in the lemmas below.

\begin{lem} \label{constr1} Let $(f_0,\alpha_0)$ be a Morse-Smale
function on a closed manifold $N$.  Let $c:N\to \RR$ be such that
$c(x)> f_0(x)$ for all $x\in N$ and let $\alpha_1$ be a second metric
on $N$.  There is a Morse-Smale function $(f,\alpha)$ on $N\times
[0,1]$ such that $f|_{N_0}=f_0$, $f|_{N_1}=c$, $\alpha$ extends the
$\alpha_i$'s and $\partial f/\partial t(x,t) >0$ for all $x\in N$,
$t>0$ $(N_i=N\times\{i\})$.  \end{lem}

\begin{proof} Let $u:[0,1]\to [0,1]$ be a $C^{\infty}$ function
such that
$$u(0)~=~1~,~u(1)~=~0~,~u'(0)~=~0~,~u'(t)<0~(t>0)~,~u''(0)<0~.$$
Consider the function
$$f(x)~=~u(t)f_0(x)+(1-u(t))c(x)~.$$
We have
$$\begin{array}{l}
\partial f/\partial t~=~u'(t)(f_0(x)-c(x))~, \\[1ex]
\partial f/\partial x~=~u(t)\partial f/\partial x +
(1-u(t))\partial c/\partial x
\end{array}$$
and this together with the fact that
$$\begin{array}{l}
(\partial^2f/\partial x^2)(x,0)~=~
(\partial^2 f_0/\partial x^2)(x)~, \\[1ex]
(\partial ^2f/\partial t^2)(x,0)~=~u''(t)(f_0(x)-c(x))
\end{array}$$
implies immediately the statement (notice that $f$ is regular on
$N_1$),
except that we also need to remark that there is a metric $\alpha$
extending the $\alpha_i$'s such that $f$ is Morse-Smale with respect to
$\alpha$.  By using a partition of unity argument it follows that there
is a metric $\alpha'$ that extends the $\alpha_i$'s.  Moreover, as
$(f_0,\alpha_0)$ is already a Morse-Smale function we can slightly
modify $\alpha'$ away from $\partial (N\times [0,1])$ to obtain
$\alpha$.
\end{proof}

\begin{lem}\label{constr2}
Consider a compact cobordism $(M;N_0,N_1)$, and suppose given
Morse-Smale functions $(f_i,\alpha_i):N_i \to \RR$ $(i=0,1)$.  Let
$g:M\to \RR$ be a Morse function with all its critical points in the
interior of $M$ and which is constant, maximal
on $N_1$, and constant, minimal on $N_0$.  For any
neighbourhood $W$ of $\partial M$ there are suitable constants $c_0$,
$c_1$ and a Morse-Smale cobordism $(f,\alpha):(M;N_0,N_1)\to \RR$
between the Morse-Smale functions $(f_0+c_0,\alpha_0)$ and
$(f_1+c_1,\alpha_1)$ such that
$$f \vert~=~g\vert~:~M\backslash W \to \RR~.$$
\end{lem}

\begin{proof} Inside the neighbourhood $W$ of $\partial M$ we can find
a tubular neighbourhood that we shall identify with $N_0\times
[0,1]\coprod N_1\times [0,1]$. We may assume that $\partial M=N_0\times
\{0\}\coprod N_1\times\{1\}$ and that $g$ is constant, regular and
equal to $k_0$ on $N_0\times \{1\}$ and is constant, regular and equal
to $k_1$ on $N_1\times \{0\}$.  Pick $c_i\in \RR$ such that
$f_0(x)+c_0<k_0$ and $f_1(x)+c_1>k_1$.  Now Lemma \ref{constr1}
provides a cobordism between $f_0+c_0$ and the constant function $k_0$
on $N_0\times [0,1]$.  By the same method as in \ref{constr1} we also
get an analogous cobordism between $k_1$ and $f_1+c_1$ on $N_1\times
[0,1]$.  Another partition of unity argument shows that these two
cobordisms can be pasted together with $g$ to provide the function $f$.
The metric $\alpha$ is obtained by the same argument as in
\ref{constr1}.
\end{proof}

\begin{rem} {\rm For the cobordism $(f,\alpha)$ constructed in
\ref{constr2} the chain complex $C(M\backslash \partial M,f)$ of
Proposition \ref{chain complex} coincides with $C(M,g)$.}\hfill\qed
\end{rem}

\begin{defi} \label{linear}
{\rm Let  $(f_0,\alpha_0)$, $(f_1,\alpha_1):N
\to \RR$ be Morse-Smale functions on a closed
manifold $N$ such that
$$f_1(x)>f_0(x)~~(x\in N)~.$$
A {\it linear} cobordism between
$(f_0,\alpha_0)$ and $(f_1,\alpha_1)$ is  a
simple cobordism (\ref{simple})
$$(f,\alpha)~:~ N \times ([0,1];\{0\},\{1\}) \to \RR$$
with $(f,\alpha)\vert_{ N\times \{i\}}=(f_i,\alpha_i)$ ($i=0,1$) such
that for all $t\in [0,1]$, $f_t$ is a convex, linear combination of
$f_0$ and $f_1$ and $\partial f/\partial t >0$ for all points $(x,t)\in
N\times (0,1)$.  \hfill\qed}
\end{defi}

\begin{lem}\label{constr3}
Any two Morse-Smale functions $(f_0,\alpha_0)$,
$(f_1,\alpha_1)$ on a closed manifold $N$ such
that
$$f_1(x)>f_0(x)~~(x\in N)$$
are related by a linear cobordism.  If $\alpha_0=\alpha_1$ there is a
simple cobordism $(f,\alpha)$ with $\alpha=\alpha_0+dt^2$ the product
metric.  If moreover $f_1=f_0+c$ with $c>0$ constant, then there exists
a linear cobordism with $\alpha$ the product metric and for any such
linear cobordism the chain map $f^{N \times[0,1]}:C(N,f_1) \to
C(N,f_0)$ is a simple isomorphism.
\end{lem}
\begin{proof} As in the proof of Lemma \ref{constr1} we consider a
$C^{\infty}$ function $v:[0,1]\to [0,1]$ such
that
$$\begin{array}{l}
v(0)=1~,~v(1)=0~,~v'(1)=v'(0)=0~,\\[1ex]
v'(t)<0~(0<t<1)~,~v''(0)<0<v''(1)~.
\end{array}$$
We define
$$f(x)~=~v(t)f_0(x)+(1-v(t))f_1$$
with the immediate consequence that, as above, one can find a metric
$\alpha$ such that the Morse-Smale function $(f,\alpha)$ satisfies the
properties required in the first part of the statement.  For the second
part, assume that $\alpha$ is the product metric $\alpha_0+dt^2$.  As
both $(f_0,\alpha_0)$, $(f_1,\alpha_1)$ are Morse-Smale we can modify
$f$ outside a neighbourhood of $\partial (N\times [0,1])$ to obtain the
desired Morse-Smale cobordism.  If $f_1=f_0+c$ this modification is not
necessary, because the function $f$ itself is already Morse-Smale with
respect to $\alpha$.  Moreover, if $x\in {\rm Crit}(f_0)$ then for all
$t\in [0,1]$ the flow induced by $-\nabla f$ is tangent to $x\times
[0,1]$.  As ${\rm ind}_f(x\times\{0\})= {\rm ind}_{f_0}(x)$,
${\rm ind}_f(x\times \{1\})= {\rm ind}_{f_0}(x)+1$ the proof is
concluded.
\end{proof}

The proof of Proposition \ref{cob_constr}
follows by applying the construction in Lemma \ref{constr3} to the
Morse-Smale functions $(f_0,\alpha_0)=(h',\gamma')$ and
$(f_1,\alpha_1)=(h,\gamma)$.
\end{proof}

\begin{rem}\label{rem:constant_shift}
{\rm  Note that if $(g,\beta)$ is a Morse-Smale function on $N$, then
the Morse
complexes $C(N,g,\beta)$ and $C(N,c(g+k),\beta)$, where $c,k\in \RR$
are constants,
are canonically identified. This implies that, on a compact manifold,
we may use
Proposition \ref{cob_constr} to compare the Morse complexes of any two
Morse
functions $h$ and $h'$ by simply adding to $h$ a sufficiently large
constant $S$
such that $h(x)+S>h'(x)$, $\forall x\in N$.}\hfill\qed
\end{rem}


\subsection{Rigidity of the Morse complex}\label{subsubsec:rigid}

We now state and prove the Rigidity Theorem \ref{theo:rig} for Morse
complexes.

Fix a Morse-Smale function $(h,\gamma)$ on a compact closed manifold
$N$. Consider
$(h',\gamma'):N \to \RR$, another Morse-Smale function on $N$, and let
$S>0$ be such
that $S>||h'-h||_{0}$. By Proposition \ref{cob_constr} and Remark
\ref{rem:constant_shift} there exists a simple Morse-Smale cobordism
$(H,\Gamma):N\times [0,1]\to \RR$ from $(h+S,\gamma)$  to
$(h',\gamma')$   such that
$$h'(x)\leq H_t(x)\leq h(x)+S~~(x \in N,t\in [0,1])$$
and the induced chain map
$$i(H)~=~H^{N\times [0,1]}~:~C(N,h,\gamma)\to  C(N,h',\gamma')$$
is a simple chain equivalence.

\newtheorem{rigtheo}[theo]{Rigidity Theorem}
\begin{rigtheo}\label{theo:rig} In the setting above
let
$$\delta~=~{\rm min}\{|h(x)-h(y)| : \ x,y\in {\rm Crit}_{\ast}(h)\ ,\
h(x)\not=h(y)\}~.$$
\begin{itemize}
\item[(i)] If $S<\delta/2$, then there exists a chain map
$j:C(N,h',\gamma')\to  C(N,h,\gamma)$ such that
$$j\circ i(H)~=~{\rm identity}~:~C(N,h,\gamma)\to  C(N,h,\gamma)~.$$
\item[(ii)] There exists a $C^2$-neighbourhood of $h$,
$\mathcal{U}_h$, such that if $h'\in
\mathcal{U}_h$ then the simple chain
equivalence $i(H)$  is a simple
isomorphism.
\end{itemize}
\hfill\qed
\end{rigtheo}

\begin{rem} {\rm (a)  The existence of a constant $S$ satisfying
condition (i) is implied by the assumption $||h-h'||_0<\delta/12$.
This means that
whenever this last condition is satisfied the Morse complex of $h$ is a
retract of
the Morse complex of $h'$.  Moreover, this relation is independent of
the metrics
used in the definition of these complexes. This also implies that the
number of
critical points of index $k$ of $h'$ is at least that of $h$.  It
should be noted
that $h$ and $h'$ do not play symmetric roles in the statement, since
$\delta$
depends on $h$.  One of the most striking features of (i) is that only
the closeness
of $h'$ to $h$ in the $C^0$ norm is required for the relation between
the Morse
complexes of $h$ and $h'$ to hold.  This allows $h'$ to have a
different number of
critical points than $h$ and, in particular, one could apply (i) to a
function $h'$
obtained from $h$ by a process analogous to the subdivision of a cell
decomposition
of a
$CW$ complex.\\
(b) Even if $f_0=f_1$ it is easy to produce examples such that the
isomorphism at (ii) is not equal to the identity.\\
(c) The condition $S< \delta/2$ appears to be optimal. It is certainly
essential
for the proof we shall present below.}
\end{rem}

In \S\ref{subsubsec:rigid} we prove \ref{theo:rig} (i), and deduce
\ref{theo:rig} (ii) as an immediate consequence.  The key new idea
appears in the proof of (i) and is as follows.
As Morse-Smale cobordisms are themselves Morse functions we may apply
the constructions in \S\ref{subsubsec:constr} to cobordisms.  This
leads in Lemma \ref{two-parameter} to a Morse function $g$ defined on
$N\times ([0,1]\times [0,1])$ having critical points only for the
second coordinate belonging to the corners of the square $[0,1]\times
[0,1]$.  The function $g$ restricts to $h$+convenient constants in the
corners $(00),(01),(11)$ and to $h'$+some constant in the corner
$(10)$.  The critical points of $g$ in the corner $(ij)$ are those of
$h$ (resp.  $h'$) but of index raised by $i+j$.  Moreover, the simple
cobordism $H$ (+ some constant) appears on the side $(10)(11)$.  On the
sides $(00)(01)$ and $(01)(11)$ we have a certain type of trivial
simple cobordism that induces the identity in terms of Morse complexes.
The purpose is to show that the restriction of $g$ on the sides
$(10)(11)$ and $(00)(10)$ induce morphisms whose composition $A$ is an
isomorphism.  For this we compare this composition with that of the
morphisms appearing on the sides $(01)(11)$ and $(00)(01)$ which, as
said before, is the identity.  The two compositions are measured by
flow lines of $g$ that join critical points of indexes differing by
$2$, that originate in $(11)$ and end in $(00)$, and are broken
precisely once at one critical point of intermediate index that belongs
to $(10)$ for the first composition and to $(01)$ for the second.  If a
flow line of $g$ originating at a critical point in $(11)$ and reaching
a critical point in $(00)$ which is broken once would break necessarily
in a point in $(10)$ or $(01)$, then the two compositions would be
equal.  This is not necessarily the case though: the breaking point can
belong also to $(00)$ or to $(11)$.  However, when the assumption on
$S$ at (i) is assumed, such a flow line with the origin in a critical
point of $x\in {\rm Crit}_i(h)\times (11)$ and ending at $y\in {\rm
Crit}_i(h)\times (00)$ with $h(y)\geq h(x)$ can only break at a point
in $(01)$ or in $(10)$ (by Lemma \ref{setting}) essentially because
along the negative-flow lines of $g$ the value of $g$ has to {\em
decrease}
and a break in a point of $(00)$ or $(11)$ would force the change in
value of $g$ to be bigger than allowed by the fact that $h(y)\geq
h(x)$.  This is enough to show that our composition $A$ is an
isomorphism even if it might not be the identity.

We now proceed to the actual proof of the Rigidity
Theorem \ref{theo:rig}.

\begin{proof} We first observe that the notion of Morse-Smale cobordism
introduced in
Definition \ref{cobordism} can be generalized in an obvious way to the
case when $N_0$ and $N_1$ have boundaries and $f_i$ is tangent to
$\partial N_i$ (in the sense that $\nabla f_i(x)\in T_{x}\partial N_i$
for all $x\in\partial N_i$, $i=0,1$).

All the previous statements have analogues in this case.  In
particular, the statements of \ref{chain complex}, \ref{chaineq},
\ref{constr1}, \ref{constr2}, \ref{constr3} remain true when assuming
that the manifold $N$ has a non-empty boundary.  A general existence
result for two-parameter Morse-Smale functions follows.

\begin{lem} \label{two-parameter}
Let $W$ be a compact manifold without boundary.  Let $(f_j,\alpha_j)$
be Morse-Smale functions on $W$, $j=0,1,2,3$ and let
$(g_{ij},\beta_{ij})$ be simple Morse-Smale cobordisms on $W\times
[0,1]$ of $(\alpha_i, f_i)$ and $(\alpha_j,f_j)$ for
$(i,j)\in\{(0,1),(0,2),(1,3),(2,3)\}$.  Assume $g_{13}>g_{02}$ and
suppose that $g_{01}$ and $g_{23}$ are linear cobordisms.  There exists
a simple Morse-Smale cobordism $(g,\beta)$ on $(W\times [0,1])\times
[0,1]$ of $(g_{02},\beta_{02})$ and $(g_{13},\beta_{13})$ which
restricts to $(g_{0+2k, 1+2k},\beta_{0+2k, 1+2k})$ on
$W\times\{k\}\times [0,1]$ $(k=0,1)$, with
$${\rm Crit}_i(g)~=~{\rm Crit}_{i-2}(f_3)\cup {\rm Crit}_{i-1}(f_2)\cup
{\rm Crit}_{i-1}(f_1)\cup {\rm Crit}_i(f_0)~.$$
\end{lem}
\begin{proof} Let
$$g_{01}~=~u(t)f_0+(1-u(t))f_1~,~g_{23}~=~r(t)f_2+(1-r(t))f_3$$
with $u,r:[0,1]\to [0,1]$ functions with the
properties described at the beginning of the
proof of Lemma \ref{constr3}. Let $w:[0,1]\to
[0,1]$ be another such function. Define
$$\begin{array}{ll}
g(x,t,\tau)~=&g_{02}(x,\tau)(u(t)w(\tau)+r(t)(1-w(\tau)) \\[1ex]
&+g_{13}(x,\tau)((1-u(t))w(\tau)+
(1-r(t))(1-w(\tau)))~.
\end{array}$$
By the same type of argument as those used before it is easy to see
that $g$ is a Morse function.  Moreover, by a partition of unity
argument one can construct a metric $\beta'$ extending the
$\beta_{ij}$'s.  A small perturbation of $\beta'$ away from the
boundary of $W\times [0,1]\times [0,1]$ leads to a new metric $\beta$
such that the pair $(\beta,g)$ satisfies the desired properties.
\end{proof}

$$\xy\xymatrix@+20pt{
 & & & \ar@{-}[r] \ar@{.}[d]_{f_1}\ar@{-}[llddd]& \ar@{-}[d]^{f_3}
\ar@{-}[llddd] \ar@{}[ld]_{g_{13}}\ar@{}[lldddd]|{g_{23}} &  \\
& &  & \ar@{.}[r]\ar@{.}[llddd] & \ar@{-}[llddd] & \\
& & & & & \\
&
\ar@{-}[r]\ar@{-}[d]_{f_0}\ar@{}[rd]^{g_{02}}\ar@{}[rruu]|{g_{01}}
&
 \ar@{-}[d]^{f_2} & & &\\
& \ar@{-}[r] & & \ar[rruuu]_t & & \\
\ar@{-}[u]^{\bf V}& \ar[r]_{\tau}& & & & }
\endxy $$

\begin{lem}\label{setting}
In the setting of Lemma \ref{two-parameter}
above, there is an $\epsilon >0$ such that for
any $u\in {\rm Crit}_{k-2}(f_3)$ and $v\in {\rm
Crit}_{k-2}(f_0)$ with $f_3(u)-f_0(v) <\epsilon$ we have
$$\sum\limits_{x\in {\rm
Crit}_{k-2}(f_1)}n^{g_{13}}(u,x)n^{g_{01}}(x,v)+
\sum\limits_{y\in {\rm
Crit}_{k-2}(f_2)}n^{g_{23}}(u,y)n^{g_{02}}(y,v)~=~0~.
\eqno{\qed}$$
\end{lem}

\begin{proof}
For $i=0,3$, let
$$\begin{array}{l}
\epsilon_i= {\rm min} \{f_i(a)-f_i(b) : a\in
{\rm {\rm
Crit}}_k(f_i),\\[2ex]
\hspace*{20mm} b\in {\rm Crit}_{k-1}(f_i), \
W^u(a)\bigcap W^s(b)\not=\emptyset, \ k\in {\bf
N} \}~.
\end{array}$$
As $W$ is compact $\epsilon_i>0$. Let $\mu>0$ be a small positive
constant
such that
$$\epsilon~=~{\rm min} \{\epsilon_0,\epsilon_3\}-\mu>0~.$$
As a consequence of the fact that $C(M,g)$ is a
chain complex we have
$$\begin{array}{l}
\sum\limits_{x\in {\rm
Crit}_{k-2}(f_1)}n^{g_{13}}(u,x)n^{g_{01}}(x,v)+
\sum\limits_{y\in {\rm
Crit}_{k-2}(f_2)}n^{g_{23}}(u,y)n^{g_{02}}(y,v)+\\[2ex]
\sum\limits_{s\in {\rm
Crit}_{k-1}(f_0)}n^g(u,s)n^{f_0}(s,v)+
\sum\limits_{l\in {\rm
Crit}_{k-3}(f_3)}n^{f_3}(u,l)n^g(l,v)~=~0~.
\end{array}$$
The condition imposed to $u$ and $v$ implies that the two last sums of
this expression vanish.  Indeed, in the sum before last the only terms
that count are those that satisfy $f_0(v)<f_0(s)<f_3(u)$.  This implies
$f_0(s)-f_0(v)<\epsilon$.  It follows that $n^{f_0}(s,v)=0$.  The
argument for the vanishing of the last sum is similar.
\end{proof}

We now use Lemma \ref{setting} to prove the Rigidity Theorem
\ref{theo:rig} (i).  We return to the setting of the statement of the
Theorem.  Thus $N$, $(h,\gamma)$, $\delta$, $S<\delta/2$ and
$(h',\gamma')$ are fixed as well as the simple Morse-Smale cobordism
$(H, \Gamma)$ of $(h',\gamma')$ and $(h+S,\gamma)$. Let $\mu>0$ such
that
$S+\mu<\delta/2$. We intend to use
Lemmas \ref{two-parameter} and \ref{setting}.  We take
$$\begin{array}{l}
W~=~N~,~(f_0,\alpha_0)~=~(h,\gamma)~,\\[1ex]
(f_1,\alpha_1)~=~(h'+S+\mu,\gamma')~,~
(f_2,\alpha_2)~=~(h+\mu,\gamma)~,\\[1ex]
(f_3,\alpha_3)~=~
(h+2S+\mu,\gamma)~,~(g_{13},\beta_{13})~=~(H+S+\mu,\Gamma)
\end{array}$$
and $(g_{02},\beta_{02})$ a linear Morse-Smale
cobordism of $(h,\gamma)$ and $(h+\mu,\gamma)$
with $\beta_{02}$ a product metric. As we have
$h'(x)+S>h(x)$ for all $x\in N$  we note that
$$g_{13}(x,\tau)~=~H(x,\tau)+S+\mu\geq
h'(x)+S+\mu>h(x)+\mu \geq g_{02}(x,\tau)~.$$ Therefore the
conditions in Lemma \ref{two-parameter} are verified.  We may also
take $(g_{23},\beta_{23})$ such that this is a linear cobordism
with $\beta_{23}$ the product metric (see Lemma \ref{constr3}).
We use Lemma \ref{two-parameter} to construct $(g,\beta)$.  Lemma
\ref{setting} can now be applied and by inspecting its proof we
see that we may take $\epsilon =\delta -\mu$.  We fix a total
order on the set ${\rm Crit}_{\ast}(h)={\rm Crit}_{\ast}(f_0)={\rm
Crit}_{\ast}(f_3)$ such that for $x,y\in {\rm Crit}_{\ast}(h)$ the
inequality $h(x)\leq h(y)$ implies $x\leq y$.  With this total
ordering we consider the matrix $A=(a_{ij})$ of the composition
$g_{01}^{N\times [0,1]}\circ g_{13}^{N\times [0,1]}$ as well as
the matrix $B=(b_{ij})$ of the composition $g_{02}^{N\times
[0,1]}\circ g_{23}^{N\times [0,1]}$.  Because both $g_{02}$ and
$g_{23}$ are linear cobordisms with the product metric it
immediately follows that $B=Id$.  We now want to observe that $A$
is an upper triangular matrix with $-1$'s on the diagonal.
Indeed, assume $x,y\in {\rm Crit}_k(h)$ correspond to elements
$i,j$ respectively in the fixed order of ${\rm Crit}_{\ast}(h)$
with $i\leq j$.  Then $f_0(x)\leq f_0(y)$ and
$$f_3(x)-f_0(y)~=~f_0(x)+2S+\mu-f_0(y)<\delta-\mu~=~\epsilon~.$$
Therefore, by Lemma \ref{setting} we have $$a_{ij}=\sum\limits_{z\in
{\rm Crit}_k(f_1)}n^{g_{13}}(x\times (1,1),z\times
(0,1))n^{g_{01}}(z\times (0,1),y\times (0,0))=-b_{ij}$$ and as
mentioned above $b_{ij}=\delta_{ij}$, the Kronecker symbol.  As a
consequence we get that $A$ has determinant equal to $\pm 1$ and is
therefore an isomorphism.  But $g_{13}^{N\times [0,1]}=H^{N\times
[0,1]}$, and this completes the proof of the Rigidity Theorem
\ref{theo:rig} (i).

We now turn to \ref{theo:rig} (ii).  If $\mathcal{U}_h$ is sufficiently
small, then
the number of critical points of $h'\in\mathcal{U}_h$ equals the number
of critical
points of $h$.  Of course, we may assume $\mathcal{U}_h$ sufficiently
small such
that $||h-h'||_0< \delta/2$ and this implies that our chain maps
$$H^{N\times [0,1]}~=~g_{13}^{N\times [0,1]}~:~ C(N,h,\gamma)\to
C(N,h',\gamma')$$
and $g_{01}^{N\times [0,1]}$ are both isomorphisms.
This completes the proof of the Rigidity Theorem
\ref{theo:rig}.
\end{proof}

\begin{rem} {\rm (a) The proof of the Rigidity Theorem \ref{theo:rig}
(i) can be used to show that for $\mathcal{U}_h$ sufficiently small
and $H'$ a second cobordism satisfying the properties of $H$ the two
isomorphisms $H^{N\times [0,1]}$ and $(H')^{N\times [0,1]}$ differ by a
nilpotent chain map.\\
(b) If $(h,\gamma)$, $(h',\gamma'):N \to \RR$
are Morse-Smale functions with $h'$
sufficiently $C^{0}$-close to $h$, the construction of
the simple cobordism $H$ in the proof of the
Rigidity Theorem \ref{theo:rig} implies that
the resulting chain map
$$i~=~H^{N\times [0,1]}~:~C(N,h,\gamma)\to  C(N,h',\gamma')$$
respects the "critical value" filtration of the chain complexes.  In
particular, if $h'=h$ and the critical points of $h$ are alone on their
critical levels, then $i$ is a simple isomorphism given in each degree
by
$$i~=~1+\hbox{\rm lower triangular matrix}~:~
C_r(N,h,\gamma)\to  C_r(N,h',\gamma')~.$$
(c) The Rigidity Theorem \ref{theo:rig} (iii) shows that a Morse
function has simple isomorphic Morse complexes for any two metrics with
respect to which it is Morse-Smale.  This can also be shown using
bifurcation
methods as in Latour \cite{Lat} (page 21).  \hfill\qed}
\end{rem}

\subsection{Mayer-Vietoris type formula for the Morse
complex}\label{gluing}

With topological cobordisms one is able to perform two natural,
geometric
operations. The first associates to a pair formed
by a cobordism $(M;N_{0},N_{1})$ and a separating closed hypersurface
$N\subset Int(M)$
the two cobordisms $(M';N_{0},N)$ and $(M'';N,N_{1})$ where
$M=M'\cup_{N}M''$.
The second is its inverse and it associates to a pair of two cobordisms
$(M';N_{0},N)$ and $(M'';N,N_{1})$ the new cobordism $(M;N_{0},N_{1})$
obtained
by pasting $M'$ to $M''$ over $N$. The fact that these two operations
are inverses
is essential for using Mayer-Vietoris type arguments to deduce
properties
of a ``long" cobordism out of knowledge about the smaller pieces in
which it can be
divided.

The purpose of this section is to put into place the same two
operations for
Morse-Smale functions defined on cobordisms and to show a
Mayer-Vietoris type theorem for Morse complexes.

Consider the pair formed by a Morse-Smale function
$(f,\alpha):(M;N_{0},N_{1}) \to \RR$ and a Morse-Smale function
$(g,\beta): N\to \RR$ with $N$ a regular hypersurface of $M$ (we
shall assume $N=f^{-1}(0)$). In this setting, the first operation
will associate to this pair two Morse-Smale functions
$h':(M';N_{0},N)\to \RR$, $h'':(M'';N,N_{1})\to \RR$ such that
$M'=f^{-1}(-\infty,0]$,$M''=f^{-1}[0,\infty)$, $h'|_{N}=kg+c'$,
$h''|_{N}=kg+c''$ with $k,c',c''$ small constants and $h'$, $h''$
equal to $f$ away from a neighbourhood of $N$ .

Conversely, given two such Morse-Smale functions $h'$ and $h''$ the
second operation
produces a Morse-Smale function $h:(M;N_{0},N_{1})\to \RR$. Obviously,
$h$ can not be
obtained by simply identifying $h'$ and $h''$ on $N$ even
if $h'|_{N}=h''|_{N}=kg+c$ (with $c=c'=c''$) because such an
identification does not
produce a Morse function.
Instead, we consider a linear Morse-Smale cobordism $L:N\times [0,1]\to
\RR$
from $kg+c+c'=L|_{N\times\{0\}}$ to $kg+c''=L|_{N\times\{1\}}$ with
$c>0$ small such that $c+c'>c''$.
We paste $L$ with $h'+c$ over $N\times\{0\}$ and then paste
the result with $h''$ at $N\times\{1\}$ and we thus  obtain a Morse-
Smale function
$h$ as desired.

Contrary to the case of usual cobordisms,
if the two operations are applied in succession starting from $f:M\to
\RR$, then the
resulting function $h$ is quite different from the initial $f$. As $h$
is in fact obtained
quite canonically from the pieces $h',h''$ in which $f$ has been
``split" we shall call such an $h$
a {\em splitting} of $f$. This terminology is also justified by the
fact that
$h$ coincides with $f$ away from a tubular
neighbourhood of $N$  and, at the same time,
any (negative) gradient flow line of $h$ that passes through a point
belonging to $M'$ (or to $M''$)
never crosses $N\times\{1/2\}$.

We formalize the construction of splittings in \S \ref{cut} and remark
that
it is immediate to express the Morse complex of $h$ purely
algebraically in terms of the complexes of $h'$ and $h''$.

It is natural to wonder why in constructing $h'$, $h''$ we do not take
simply $h'=f'$, $h''=f''$
where $f'=|_{M'}$, $f''=f|_{M''}$. The reason is that the simplest form
of a Mayer-Vietoris type
formula for Morse complexes should provide
the Morse complex of $f$ out of the Morse complexes of the
pieces $h'$ and $h''$. However, the Morse complexes of $f'$ and $f''$
are inappropriate for this
task because they do not encode the flow lines of $-\nabla f$ that
cross from $M''$ into $M'$.
On the other hand, the gluing
Theorem \ref{theo:glueadapt} in \S\ref{assembly} shows that if $g$
satisfies a
certain, generic, technical condition, then the Morse complexes of $h'$
and $h''$
{\em can} be ``glued" together in a canonical way to give the Morse
complex of $f$. The isomorphism
type of the Morse complex of $f$ can be recovered out of any splitting
of $f$
even if the technical condition is not satisfied.

\subsubsection{Splittings}\label{cut}

In this section we consider a Morse-Smale function
$$(f,\alpha):(M;N_0,N_1) \to \RR$$
with $f(x_0)<0$ for all $x_0 \in N_0$, $f(x_1)>0$ for all $x_1 \in
N_1$,
and $0 \in \RR$ a regular value of $f$. We write
$$N~=~f^{-1}(0) \subset M~.$$
The restrictions
$$\begin{array}{l}
(f',\alpha')~=~(f,\alpha)\vert~:~(M';N_0,N)~=~f^{-1}(-\infty,0]
\to
\RR~,\\[1ex]
(f'',\alpha'')~=~(f,\alpha)\vert~:~(M'';N,N_1)~=~f^{-1}[0,\infty)\to
\RR
\end{array}$$
are Morse-Smale functions with
$$(f,\alpha)~=~(f',\alpha') \cup (f'',\alpha''):
(M;N_0,N_1)~=~(M';N_0,N)\cup_N (M'';N,N_1) \to
\RR~.$$

\begin{defi} \label{gluchain}
{\rm The {\it attaching chain map}
$$\phi~:~C(M'',f'',\alpha'')_{*+1} \to C(M',f',\alpha')$$
is the chain map given by
$$\phi~:~C(M'',f'',\alpha'')_{i+1} \to C(M',f',\alpha')_i~;~
x \mapsto \sum\limits_{y \in {\rm
Crit}_i(f')}n^{f,\alpha}(x,y)y~.\eqno{\qed}$$}
\end{defi}

\begin{prop} \label{glu1}
The Morse complexes fit into a short exact
sequence
$$0 \to C(M',f',\alpha') \to C(M,f,\alpha) \to C(M'',f'',\alpha'') \to
0$$ with
$$C(M,f,\alpha)~=~C(\phi)$$
the algebraic mapping cone of the attaching
chain map
$$\phi~:~C(M'',f'',\alpha'')_{*+1} \to C(M',f',\alpha')~.$$
\end{prop}
\begin{proof} Immediate from the definitions.
\end{proof}

The key technical construction of this section is contained in the
following definition.  To fix ideas we will proceed from here under the
following assumption: there exists a neighbourhood
$N\times[-\delta,\delta]\subset M$ of $N$ such that on this
neighbourhood $f$ has the form $f(x,t)=t$ and the metric $\alpha$ has
the form $dt^2+\beta $, $\beta=\alpha|_N$.  The existence of such a
parametrization of $f$ around $N$ is clear from the fact that $N$ is a
regular level hypersurface of $f$.  Moreover, because
$\partial/\partial t$ is orthogonal to $N\times \{t\}$ the metric
$\alpha$ has the form $u(x,t)dt^2+v(x,t)\beta$ (with $u,v$ smooth
positive functions) on our neighbourhood.  This means that by replacing
$\alpha$ in this neighbourhood with the product metric we do not
perturb
the flow lines of $f$.  Therefore, our assumption is in no way
restrictive.

\begin{defi} \label{splitting} {\rm An $(\epsilon,\tau)$-{\it
splitting} of a Morse-Smale function
$$(f,\alpha):(M;N_0,N_1) \to \RR$$
along a Morse-Smale function $(g,\beta):N=f^{-1}(0)\to \RR$
($\beta=\alpha\vert)$ is a Morse-Smale function $(h,\alpha):M \to \RR$
such that for $(h',\alpha')=(h,\alpha)\vert_{M'}$, $(h'',\alpha'')
=(h,\alpha)\vert_{M''}$ we have
\begin{itemize}
\item[(i)] $h'\vert_N=h''\vert_N=\tau g:N \to \RR$ with $\tau \in
\RR^{+}$.
\item[(ii)] $(h,\alpha)=(f,\alpha)$ except in a small tubular
neighbourhood
$W=N \times [-\epsilon,\epsilon] \subset M$,
$\epsilon\in\RR^{+}$, of $N=N \times \{0\}
\subset M$.
\item[(iii)] With respect to this parametrization $f(x,t)=t$ for all
$(x,t)\in N\times
[-\epsilon,\epsilon]$ and
$$\begin{array}{l}
{\rm Crit}_i(h)~=~{\rm Crit}_i(f)\cup ({\rm
Crit}_{i-1}(g)\times \{-\epsilon/2\}) \cup
({\rm Crit}_i(g)\times
\{\epsilon/2\})~,\\[1ex]
{\rm Crit}_i(h')~=~{\rm Crit}_i(f')\cup
({\rm Crit}_{i-1}(g)\times \{-\epsilon/2\})~,\\[1ex]
{\rm Crit}_i(h'')~=~{\rm Crit}_i(f'')\cup({\rm
Crit}_i(g)\times \{\epsilon/2\})~.
\end{array}$$
\item[(iv)] The restrictions of $(h,\alpha)$ to the submanifolds
$$\begin{array}{l}
W_{\epsilon}~=~N \times [-\epsilon/2,\epsilon/2] \subset W~,\\[1ex]
M'_{\epsilon}~=~M' \backslash \big(N \times
(-\epsilon/2,0]\big) \subset
M'~,\\[1ex]
M''_{\epsilon}~=~M'' \backslash \big(N \times
[0,\epsilon/2)\big) \subset M''
\end{array}$$
{\rm (}which are copies of $W,M',M''$
respectively{\rm )} are Morse-Smale cobordisms
$$\begin{array}{l}
(g_{\epsilon},\beta_{\epsilon})~=~(h,\alpha)\vert~:~
(W_{\epsilon};N\times \{\epsilon/2\},N\times
\{-\epsilon/2\}) \to
\RR~,\\[1ex]
(h'_{\epsilon},\alpha'_{\epsilon})~=~(h',\alpha')\vert~:~
(M'_{\epsilon};N_0,N\times \{-\epsilon/2\}) \to \RR~,\\[1ex]
(h''_{\epsilon},\alpha''_{\epsilon})~=~(h'',\alpha'')\vert~:~
(M''_{\epsilon};N\times \{\epsilon/2\},N_1) \to
\RR~.
\end{array}$$
\item[(v)] The cobordism $(g_{\epsilon},\beta_{\epsilon})$ is linear.
\qed
\end{itemize}}
\end{defi}

\begin{figure}
\centering
\includegraphics[width=4.7in]{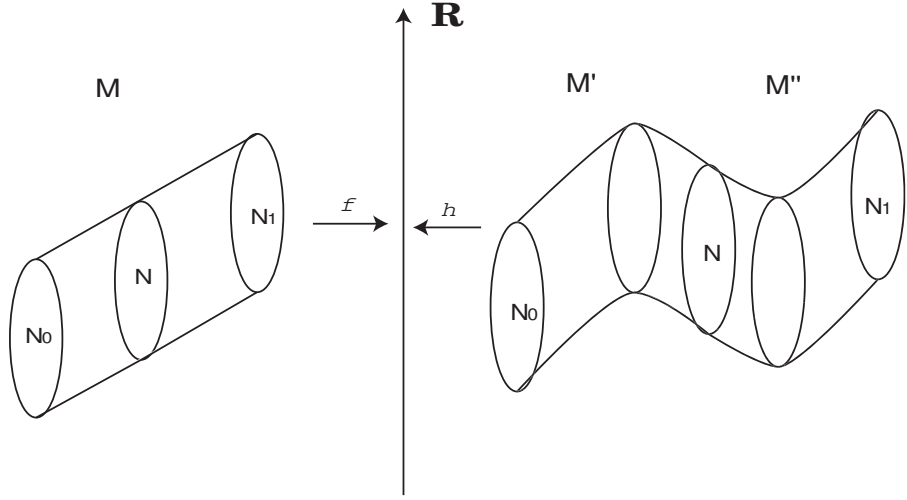}
\caption{ $f$ and its splitting $h$.} \label{fig:twodisks}
\end{figure}

\

The key properties of a splitting are collected next.

\begin{prop} \label{theta1}
Let $(f,\alpha):(M;N_0,N_1) \to \RR$
be a Morse-Smale function, with $0 \in \RR$ a regular value.\\
{\rm (i)} For every Morse-Smale function
$(g,\beta):N=f^{-1}(0) \to \RR$ with
$\beta=\alpha\vert_N$ and every $\delta >0$
there exists an $(\epsilon, \tau)$-splitting
$(h,\alpha)$ of $(f,\alpha)$
along $(g,\beta)$ with $||h-f||_0\leq\delta$.\\
{\rm (ii)} The Morse complexes of a splitting
of $(f,\alpha)$ along $(g,\beta)$
$$\begin{array}{l}
(h,\alpha)~=~(h'_{\epsilon},\alpha'_{\epsilon})
\cup (g_{\epsilon},\beta_{\epsilon})
\cup(h''_{\epsilon},\alpha''_{\epsilon})~:\\[2ex]
(M;N_0,N_1)~=~ (M'_{\epsilon};N_0,N\times
\{-\epsilon/2\}) \cup
(W_{\epsilon};N\times \{-\epsilon/2\},N\times \{\epsilon/2\})\\[2ex]
\hskip200pt  \cup (M''_{\epsilon};N\times
\{\epsilon/2\},N_1) \to \RR~,
\end{array}$$
are given by
$$\begin{array}{l}
C(M,h,\alpha)~=~C(\theta)~,\\[1ex]
C(M'_{\epsilon},h'_{\epsilon},\alpha'_{\epsilon})~=~C(M',h',\alpha')~
=~C(\theta')~,\\[1ex]
C(M''_{\epsilon},h''_{\epsilon},\alpha''_{\epsilon})~=~C(M'',h'',
\alpha'')~=~C(\theta'') \\[1ex]
C(W_{\epsilon},g_{\epsilon},\beta_{\epsilon})~=~
C(1:C(N,g,\beta) \to C(N,g,\beta))~
\end{array}$$
with
$$\begin{array}{l}
\theta~:~C(M',h',\alpha')_{*+1} \to C(M'',h'',\alpha'')~,\\[1ex]
\theta'~:~C(N,g,\beta)\to C(M',f',\alpha')~,\\[1ex]
\theta''~:~C(M'',f'',\alpha'')_{*+1}  \to C(N,g,\beta)
\end{array}$$
the chain maps {\rm (}up to sign{\rm )} defined by
$$\begin{array}{l}
\theta~:~C(M',h',\alpha')_{i+1}~=~
C(M',f',\alpha')_{i+1}\oplus C(N,g,\beta)_i \\[1ex]
\hspace*{15mm}\to
C(M'',h'',\alpha'')_i~=~C(N,g,\beta)_i\oplus
C(M'',f'',\alpha'')_i ~;\\[1ex]
\hspace*{15mm}(x,z\times \{-\epsilon/2\})
\mapsto (z\times \{\epsilon/2\},0)~~
(x \in {\rm Crit}_{i+1}(f'),z \in {\rm Crit}_i(g))~,\\[2ex]
\theta'~:~C(N,g,\beta)_i \to
C(M',f',\alpha')_i~;~ z \mapsto \sum\limits_{y
\in {\rm Crit}_i(f')}
n^{h',\alpha'}(z\times \{-\epsilon/2\},y)y~,\\[3ex]
\theta''~:~C(M'',f'',\alpha'')_{i+1} \to
C(N,g,\beta)_i~;~ x \mapsto \sum\limits_{z \in
{\rm Crit}_i(g)}
n^{h'',\alpha''}(x,z\times
\{\epsilon/2\})z~.
\end{array}$$
In particular, there are defined exact sequences
$$\begin{array}{l}
0 \to C(M'',h'',\alpha'') \to C(M,h,\alpha) \to C(M',h',\alpha') \to
0~,\\[1ex]
0 \to C(M',f',\alpha') \to C(M',h',\alpha') \to C(N,g,\beta)_{*-1} \to
0~,\\[1ex]
0 \to C(N,g,\beta) \to C(M'',h'',\alpha'') \to C(M'',f'',\alpha'') \to
0~.
\end{array}$$
{\rm (iii)}
Let $p_h:C(W_{\epsilon},g_{\epsilon},\beta_{\epsilon}) \to
C(M,h,\alpha)$
be the chain map defined by the split injections
$$\begin{array}{l}
p_h~:~C(W_{\epsilon},g_{\epsilon},\beta_{\epsilon})_i~
=~C(N,g,\beta)_{i-1} \oplus C(N,g,\beta)_i\to \\[2ex]
C(M,h,\alpha)_i\\[2ex]
=~C(M',f',\alpha')_i \oplus
C(N,g,\beta)_{i-1} \oplus C(N,g,\beta)_i \oplus
C(M'',f'',\alpha'')_i~;\\[2ex]
\hspace*{50mm} (x,y) \mapsto (\theta'(y),x,y,0)~.
\end{array}$$
The  algebraic mapping cone of $p_h$ is the Morse complex of a
splitting
$(h,\alpha):M \to \RR$
$$C(p_h)~=~C(M,h,\alpha)~.$$
The cokernel of $p_h$ is the algebraic mapping cone of
$\theta'\theta'':C(M'',f'',\alpha'')_{*+1} \to C(M',f',\alpha')$
$$0 \to C(W_{\epsilon},g_{\epsilon},\beta_{\epsilon})
\xrightarrow[]{p_h} C(M,h,\alpha) \xrightarrow[]{j_h}
C(\theta'\theta'') \to 0~,$$
with $j_h$ the natural projection.
\end{prop}
\begin{proof}
(i) This follows easily from the constructions in
\S\ref{subsubsec:constr}. We first choose $a>0$ such that $[-a,a]$
contains only regular values of $f$, with
$$f\vert~:~W~=~f^{-1}([-a,a])~=~N \times [-a,a] \to \RR~;~(x,t) \mapsto
t$$
and the assumption that $\alpha$ is the product metric here.  After
possibly multiplying $g$ with a sufficiently small constant we may
assume that the image of $g$ belongs to $[-b,b]$ with $b\in (0,1)$ such
that $b\ll a$.  Let
$$\begin{array}{l}
N_{\delta}~=~N\times \{\delta\}\subset W~=~N\times [-a,a]~,\\[1ex]
W(u,v)~=~N\times [u,v]\subset W~.
\end{array}$$
By using Lemmas \ref{constr1}, \ref{constr2} for each $0<\epsilon < b$
we construct a Morse-Smale function $f_{\epsilon}:M\to \RR$ of $f$
along $g$ such that : $f_{\epsilon}$ agrees with $f$ outside
$W(-\epsilon,\epsilon)$; the gradient of $f_{\epsilon}$ is everywhere
tangent to $N_{\epsilon /2}$ and $N_{-\epsilon /2}$,
$$f_{\epsilon}|_{N_{\epsilon /2}}~=~\epsilon(g+k)~,~
f_{\epsilon}|_{N_{-\epsilon /2}}~=~\epsilon(g+k')$$
with $k'>k$ small; if $x\in {\rm Crit}_{f_{\epsilon}}\bigcap
N_{\epsilon /2}$ then ${\rm ind}_{f_{\epsilon}}(x)={\rm ind}_g(x)$ and
if $x\in {\rm Crit}_{f_{\epsilon}}\bigcap N_{-\epsilon /2}$, then ${\rm
ind}_{f_{\epsilon}}(x)={\rm ind}_g(x)+1$; on the cobordism
$(W(-\epsilon/2,\epsilon/2);\allowbreak N_{-\epsilon /2},N_{\epsilon
/2})$ the function $f_{\epsilon}$ restricts to a linear cobordism.
The construction of the function $f_{\epsilon}$ is done by applying
Lemma
\ref{constr1} to the two cobordisms $W(-\epsilon,-\epsilon /2)$ and
$W(\epsilon/2, \epsilon)$ and Lemma \ref{constr2} to the cobordism
$W(-\epsilon/2, \epsilon/2)$.  In particular, $f_{\epsilon}$ is a
linear combination of the form $u(t)(-\epsilon)+(1-u(t))\epsilon(g+k')$
inside $N\times [-\epsilon,-\epsilon/2]$ with
$u:[-\epsilon,-\epsilon/2]\to \RR$ an appropriate function with
$u(-\epsilon)=1$ and $u(-\epsilon/2)=0=u'(-\epsilon/2)$, $u''(-
\epsilon/2)>0$
(see \ref{constr1}); a similar
linear combination is valid inside $N\times [\epsilon/2,\epsilon/]$.
To understand the pasting between the cobordism $W(-
\epsilon/2,\epsilon/2)$
and, for example, $W(-\epsilon, -\epsilon/2)$ notice that on $W(-
\epsilon/2,\epsilon/2)$
we have a cobordism from $(\epsilon(g+k'),N_{-\epsilon/2})$ to
$(\epsilon(g+k),N_{\epsilon/2})$. In other words, on both sides of
$N_{-\epsilon/2}$
we have cobordisms having $N_{-\epsilon/2}$ as ``high end". This allows
the pasting
to occur by simply using in the formula above $u$ defined on a larger
interval, say $u:[-\epsilon,0]\to \RR$, and
using the same formula to define $f_{\epsilon}$ on both sides of $N_{-
\epsilon/2}$.
For a generic choice of $g$ the function $f_{\epsilon}$ is already
Morse-Smale with respect to the metric $\alpha$.  In general, to obtain
a function $h$ such that $(h,\alpha)$ is Morse-Smale we might need to
still slightly perturb $f_{\epsilon}$ inside $N\times
([-\epsilon,-3\epsilon/4]\cup [3\epsilon/4,\epsilon])$.  We may assume
that the construction has been made such that $||f-f_{\epsilon}||_0\leq
5\epsilon$ and by adjusting the various relevant constants we see that
$(h,\alpha)$ is a splitting of $f$ along $g$ with the required
properties.\\
(ii)+(iii) Immediate from the definition and construction of a
splitting.
\end{proof}

\subsubsection{Assembly of the Morse complex}\label{assembly}

Given the algebraic data extracted from a splitting of a Morse-Smale
function $(f,\alpha):(M;N_0,N_1) \to \RR$ we now describe how to
reconstruct the Morse complex of $f$ from the Morse complexes of $h'$
and $h''$.  One expects the composite chain map
$$\theta'\theta''~:~C(M'',f'',\alpha'')_{*+1} \to C(M',f',\alpha')$$
to be a good approximation of the attaching chain map $\phi$.

\begin{defi} \label{adapted}
{\rm A Morse-Smale function $(g,\beta):N=f^{-1}(0) \to \RR$
is {\it adapted} to splitting
$(f,\alpha):(M;N_0,N_1) \to \RR$ at $N \subset
M$ if
\begin{itemize}
\item[(i)] $(g,\beta):N \to \RR$ is Morse-Smale, with
$\beta=\alpha\vert_N$.
\item[(ii)] For any $p\in {\rm Crit}_i(f)$
$$W_f^u(p)\bigcap N\subset
\bigcup_{z\in {\rm Crit}_{\leq
i-1}(g)}W_g^u(z)~.$$
\item[(iii)] For any $q\in {\rm Crit}_{\ast}(f)$ and any $z\in {\rm
Crit}_{\ast}(g)$ the intersections \\
$W^u_f(q)\bigcap W^s_g(z)$ and $W^u_g(z)\bigcap
W^s_f(q)$ are transversal.
\end{itemize}
\hfill\qed}
\end{defi}

\newtheorem{glutheo}[theo]{gluing Theorem}

\begin{glutheo} \label{theo:glueadapt}
Let $(f,\alpha):(M;N_0,N_1) \to \RR$ be a Morse-Smale function, with
$N=f^{-1}(0)$ a regular hypersurface.\\
{\rm (i)} Let $(g,\beta):N\to \RR$ be a Morse-Smale function.  For any
$(\epsilon,\tau)$-splitting $(h,\alpha)$ of $f$ along $g$ with
$\epsilon,\tau$ sufficiently small there exists a
chain homotopy $\psi:\phi \simeq \theta'\theta''$,
defining a simple isomorphism
$$\begin{pmatrix} 1 & \psi \\ 0 & 1 \end{pmatrix}~:~
C(M,f,\alpha)\xrightarrow[]{\cong} {\rm coker}(p_h)~.$$
{\rm (ii)}  There exist Morse-Smale functions
$(g,\beta)$ adapted to splitting $(f,\alpha)$.\\
{\rm (iii)} If $(g,\beta)$ is adapted to
splitting $(f,\alpha)$ at $N$, then there exist splittings
$h$ of $f$ along $g$ with $\psi=0$.
Equivalently, we have the equality $\phi=\theta'\theta''$ or,
explicitly:
$$n^{f,\alpha}(x,y)~=~\sum\limits_{z\in {\rm Crit}_{i-1}(g)}
n^{h'',\alpha''}(x,z\times\{\epsilon/2\})
n^{h',\alpha'}(z\times\{-\epsilon/2\},y)$$
for each pair $x\in {\rm Crit}_i(f'')$, $y\in
{\rm Crit}_{i-1}(f')$.\hfill$\qed$
\end{glutheo}

\begin{rem}\label{rem:glue1} {\rm
To our knowledge, a gluing formula as the one above has not yet
appeared in the literature.  However, a closely related result is that
of Laudenbach \cite{Lau} and has been extended to $S^{1}$-valued
functions by Hutchings \cite{H}.  To see the relation with these
results notice that, in the situation of Proposition \ref{theta1}, the
based f.g.  free $\ZZ[\pi]$-module chain complex
$${\rm coker}(p_h)~=~C(\theta'\theta'')$$ is the Morse complex
$C(M,\widehat{f},\alpha)$ of a Morse-Smale function
$(\widehat{f},\alpha):M \to \RR$ obtained from $(h,\alpha)$ by
cancelling
the pairs of critical points added in the construction of $h$ from $f$,
with $${\rm Crit}_*(\widehat{f})~=~{\rm Crit}_*(f)~.$$
Indeed, the results of Laudenbach describe the modifications occurring
in the Morse complex after the birth (or death) of two mutually
annihilating critical points of successive indexes.  By applying
iteratively his result when cancelling all such pairs of critical
points of $h$ that belong to $W_{\epsilon}$ one obtains a function
$\widehat{f}$ whose Morse complex does equal $C(\theta'\theta'')$.
Moreover, this function may be assumed to be $C^{0}$-close to $f$.
Therefore, this argument, together with the rigidity theorem, is
sufficient to prove the point (i) of the theorem.  On the other hand,
this is not enough to show the point (iii).  The reason is that
``cancelling of critical points" is a non-unique operation and there is
no way to insure in general that the function $\widehat{f}$ obtained at
the end of this process is $C^{2}$-close to $f$ and has therefore a
Morse
complex identical to that of $f$.
Because of this we shall prove the whole theorem without the use of
bifurcation arguments (as we shall see, the point (i) is in fact
immediate by using Morse cobordisms).

Other constructions related to the gluing formula appear in
Pajitnov \cite{P}.  The role of our formula is played there by an
analysis of the intersections of stable and unstable manifolds of $f$
with the stable and unstable manifolds of $g$.\hfill\qed}
\end{rem}

\begin{proof}
(i) This point is quite immediate.  Assume $(h,\alpha)$ is an
$(\epsilon,\tau)$-splitting of $f$ along $g$ as in the statement.  By
the constructions in \S\ref{subsubsec:constr} there is a Morse
cobordism $H$ between $(h,\alpha)$ and $(f+c,\alpha)$.  If $\epsilon$
is small enough, then $h$ is sufficiently $C^0$ close to $f$ such that
$$\begin{array}{l}
{\rm min}\{f({\rm Crit}_{\ast}(f''))\}\geq {\rm
max}\{h({\rm Crit}_{\ast}(h)-{\rm Crit}_{\ast}(f))\}~,\\[1ex]
{\rm min}\{h({\rm Crit}_{\ast}(h)-{\rm Crit}_{\ast}(f))\}>c+ {\rm
max}\{f({\rm Crit}_{\ast}(f'))\} \end{array}$$
 and we may construct $H$ such that it restricts to linear cobordisms
with a product metric on $M'\backslash (N\times (-\epsilon,\epsilon))$
and on $M''\backslash (N\times (-\epsilon,\epsilon))$ (we use the
notations in Proposition \ref{theta1}).  Denote by $G:C(M,f,\alpha)\to
C(M,h,\alpha)$ the chain map induced by this cobordism. Due to the
linearity of $H$ on the complement of $N\times (-\epsilon,\epsilon)$
the composite
$$j'~=~j_h\circ G~:~C(M,f,\alpha) \to {\rm coker}(p_h)$$
is a simple isomorphism of the form
$$\begin{array}{l}
j'~=~\begin{pmatrix} 1 & \psi \\ 0 & 1 \end{pmatrix}~:~
C(M,f,\alpha)_i~=~C(M',f',\alpha')_i \oplus
C(M'',f'',\alpha'')_i\\[1ex]
\hskip120pt \to {\rm coker}(p_h)_i~=~C(M',f',\alpha')_i \oplus
C(M'',f'',\alpha'')_i
\end{array}$$
for a chain homotopy
$$\psi~:~\phi~ \simeq~ \theta'\theta''~:~C(M'',f'',\alpha'')_{*+1}
\to C(M',f',\alpha')~.$$
(ii) Clearly, a Morse-Smale function
$(g,\beta):N \to \RR$ is not in general adapted to splitting
$(f,\alpha):M \to \RR$.
We obtain an adapted $(g,\beta)$ by starting with an arbitrary
$(g',\beta):N \to
\RR$ and adding pairs of mutually cancelling critical points such that
the unstable
manifolds of $g$ give a $CW$ decomposition of $N$ which is sufficiently
fine to
ensure that Definition \ref{adapted} (ii) is verified.
In particular, this means that $W^{u}_{f}(R)\bigcap
W^{u}_{g'}(R')=\emptyset$
whenever $ind_{g'}(R')\geq ind_{f}(R)$. After further subdivision we
may also assume that
for any pair $R'\in {\rm Crit}_{\ast}(g')$, $R\in {\rm Crit}_{\ast}(f)$
we either
have $W^{u}_{f}(R)\cap
W^{u}_{g'}(R')=\emptyset$, or $W^{u}_{g'}(R')\subset W^{u}_{f}(R)$.
We then modify $g'$ such as to
perturb the stable manifolds $W^{s}_{g'}(R'')$  to render them
transverse to $W^{u}_{f}(R)$ for all pairs $R''\in {\rm
Crit}_{\ast}(g'),
R\in {\rm Crit}_{\ast}(f)$.
This modification is performed successively for $R''$ of increasing
indexes and it
may be achieved without modifying those unstable manifolds
$W^{u}_{g'}(R')$ which intersect some
unstable manifold of $f$. This happens because if $W^{u}_{g'}(R')$
intersects
$W^{u}_{f}(R)$, then $W^{u}_{g'}(R')\subset W^{u}_{f}(R)$ and as $g'$
is Morse-Smale
we also have that $W^{s}_{g'}(R'')$ is transversal to $W^{u}_{g'}(R')$
(for all
$R',R''\in {\rm Crit}_{\ast}(g')$). This means that if
$W^{s}_{g'}(R'')$ intersects $W^{u}_{f}(R)$
in a point which belongs also to $W^{u}_{g'}(R')$, then
$W^{s}_{g'}(R'')$ is already
transversal to $W^{u}_{f}(R)$ at this point.
We also need to modify our function $g'$ such as to obtain one which
also
has the property that the second type of intersections in
Definition \ref{adapted} (iii) are transversal. We perturb
successively  each unstable manifold $W^{u}_{g'}(R')$, $R'\in{\rm
Crit}_{\ast}(g')$
such as to make $W^{u}_{g'}(R')$ transversal to the stable manifolds of
$f$.
Obviously, we need to insure that condition (ii) in Definition
\ref{adapted} is preserved and therefore the critical points  $R'$ such
that
$W^{u}_{g'}(R')$ intersects some $W^{u}_{f}(R)$ need to be discussed
separately.
In this case we have $W^{u}_{g'}(R')\subset
W^{u}_{f}(R)$. Suppose that $ind_{f}(R)=ind_{g'}(R')+1$, then as
$W^{u}_{f}(R)$
is transverse to all stable manifolds $W^{s}_{f}(Q)$ (by the Morse-
Smale condition) it also follows
that $W^{u}_{g'}(R')$ is transverse to these stable manifolds and no
modification of
$W^{u}_{g'}(R')$ is needed. In case $ind_{g'}(R')< ind_{f}(R)+1$ it
follows from condition
(ii) in Definition \ref{adapted} that
$W^{u}_{g'}(R')\subset W^{u}_{f}(R)\cap N\subset \bigcup_{ind_{g'}(z)=
ind_{f}(P)-1}\overline{W^{u}_{g'}(z)}=T$. We then perturb
$W^{u}_{g'}(R')$ inside
$W^{u}_{f}(R)$ such as to obtain the needed transversality. This will
also perturb
the higher dimensional $W^{u}_{g'}(z)$'s that appear in the union $T$
above but can be made
such that $T$ does not diminish and thus (ii) of Definition
\ref{adapted} continues to be
satisfied.

It is useful to note that one can obtain in this way an adapted $g$
that
is arbitrarily close in $C^0$ norm to the function
$g'$ that was given initially.\\
(iii) We now assume that $(g,\beta)$ is adapted to splitting
$(f,\alpha)$.  Recall the functions $f_{\epsilon}$ constructed in the
proof of Proposition \ref{theta1}.  We shall work below with functions
of this form and with the notations of that proof.  In particular,
$$W~=~f^{-1}([-a,a])~,~W(u,v)~=~N\times [u,v]~=~f^{-1}([u,v])~,~
N~=~f^{-1}(0)~.$$
It is immediate to see that because of the transversality properties of
the adapted $g$, the function $f_{\epsilon}$ is Morse-Smale with
respect to $\alpha$.  Therefore $(f_{\epsilon},\alpha)$ is itself a
splitting and we will show below that any such splitting with
$(\epsilon,\tau)$ sufficiently small verifies the conclusion.\\
\indent
By inspecting the construction in \ref{theta1}, we see that, in each
point of $W(-\epsilon,\epsilon)$, the tangent vector $\nabla
f_{\epsilon}$ projected onto $N$ has the same direction as $\nabla g$.
We intend to compare the Morse complex of $f$ with that of
$f_{\epsilon}$ using a special Morse cobordism between $f_{\epsilon}$
and $f+c$.  For this we start with a linear cobordism $F:M\times
[0,1]\to \RR$ of Morse-Smale functions between $f$ and $f+c$ with $c>0$
fixed.  We assume $F_0=f$ and $F_1=f+c$.  We denote by $\tau\in [0,1]$
the deformation parameter of $F$.  We may also assume that in a
neighbourhood $M\times [0,1/3]$ of $M\times \{0\}$ the function $F$ has
the form $F=f+\tau^2$.  We shall identify below $M$ and $M\times \{0\}$
(for example we write $F|_{M\times \{0\}}=f$ etc).\\
\indent We want to show that we can modify $F$ only inside
$M^{\epsilon}=W(-\epsilon,\epsilon)\times [0,(5\epsilon)^{1/3}]$
such that the resulting function $F^{\epsilon}$ will be a
Morse-Smale cobordism of $f_{\epsilon}$ and $f+c$.  By a similar
method to that used in \ref{constr1} we can define $F^{\epsilon}$
inside $M^{\epsilon}$ by
$$F^{\epsilon}|_{M^{\epsilon}}~=~\tau^2+v(\tau)f+(1-
v(\tau))f_{\epsilon}$$
where
\begin{itemize}
\item[] $v:[0,1]\to [0,1]$ is $C^{\infty}$;
\item[] $v$ is constant equal to $1$ if $\tau> (5\epsilon)^{1/3}$;
\item[] $v(0)=v'(0)=v''(0)=0$;
\item[] $|v'(\tau)5\epsilon | < 2\tau$ for all $\tau$.
\end{itemize}
(Such a function exists, it can be modelled on $\tau^3/5\epsilon$) and
$F^{\epsilon}=F$ outside of $M^{\epsilon}$.  The function
$F^{\epsilon}$ is smooth because $f^{\epsilon}=f$ outside of
$W(-\epsilon,\epsilon)$ and it is clear that it is Morse.  Note for
each point of $M^{\epsilon}$ the gradient of $F^{\epsilon}$ with
respect to the product metric $\alpha + d\tau ^2$ decomposes as an
orthogonal sum of a component tangent to $N_0$ which is a multiple of
the
gradient of $g$ and two other components, one in the direction of
$\partial /\partial t$, and the other into that of $\partial /\partial
\tau$.  We can change slightly the metric $\alpha + d\tau ^2$ inside a
set $M\times [l,l']$ with $(5\epsilon)^{1/3} < l < l'< 1$ thus getting
a metric $\alpha_{\epsilon}$ sufficiently close to $\alpha$ such that
with respect to this new metric $F$ remains Morse-Smale and its complex
is not modified and $F^{\epsilon}$ becomes also Morse-Smale.
In fact, the set of metrics $\alpha^{\ast}$ such that $F$ remains
Morse-Smale
with respect to $\alpha^{\ast}$ and has the same Morse complex as with
respect
to $\alpha$ is dense and open in some neighborhood of $\alpha$.
Therefore,
if we now consider the sequence of functions $F^{1/n}$, $n\in {\bf
N}^{\ast}$, $n\geq n_0\gg 1$, we see that there is a metric $\alpha'$
on $M\times [0,1]$ such that $\alpha'$ equals $\alpha + d\tau ^2$
outside $M\times [1/3,2/3]$, for all $n$, $F^{1/n}$ is Morse-Smale with
respect to $\alpha'$ and $F$ is also Morse-Smale with respect to
$\alpha'$ and has the same Morse complex as it has with respect to
$\alpha + d\tau^2$ (in fact, the set of such metrics, as a countable
intersection of open
dense sets, is even dense in the mentioned neighborhood
of $\alpha$). It will be convenient to use the convention that
$n^t(p,q)=0$ whenever the difference of indexes between the critical
points $p$ and $q$ of a Morse-Smale function $t$ is strictly greater
than $1$, or if one of $p,q$ is not a critical point of $t$.

\begin{lem}\label{deform} In the setting above there is some
$n\in {\bf N}$ such that in the Morse complex
of $F^{1/n}$ we have:
\begin{itemize}
\item[(i)] $n^{F^{1/n}}(p,q)=n^F(p,q)$ if $p,q\in {\rm
Crit}(F^{1/n})\bigcap {\rm Crit}(F)$, $p\in
M\times \{1\}$.
\item[(ii)] $n^{F^{1/n}}(p,q)=0$ for $p,q\in {\rm Crit}(F^{1/n})$,
 with $p\in M\times \{1\}$ and $q\in
N_{1/(2n)}$.
\item[(iii)] $n^{F^{1/n}}(p,q)=0$ for $p,q\in {\rm Crit}(F^{1/n})$
 with $p\in M\times \{1\}$,\\ $q\in
N_{-1/(2n)}\bigcup N_{1/(2n)}$ and
$F^{1/n}(p)<0$.
\end{itemize}
\end{lem}

Assuming this Lemma, here is the end of the proof of the gluing
Theorem
\ref{theo:glueadapt}.  To simplify notation we denote $G=F^{1/n}$,
$h=f_{1/n}$ with
$n$ satisfying the conclusion of \ref{deform}.  Denote by
$\overline{F},
\overline{G}$ respectively the chain maps $F^{M\times [0,1]}$ and
$G^{M\times
[0,1]}$.  If $x\in M\times\{0\}\subset {\rm Crit}(f)\bigcap M\times
[0,1]$ we let
$x'=x\times \{1\}$.  As $F$ is a linear cobordism we have
$\overline{F}(x')=x$.
Consider now the function $h$.  It restricts to a linear cobordism on
$W(-1/(2n),1/(2n))$.  If $y\in {\rm Crit}(h)\bigcap N_{-1/(2n)}$, then
let $\theta
y=h^{W(-1/(2n),1/(2n))}(y)$.  Clearly,
$$\theta y ~=~\overline{u}\times \{1/(2n)\}~~,~~
y~=~\overline{u}\times \{-1/(2n)\}$$
with $\overline{u}$ a critical point of $g$.  Assume $x,y\in {\rm
Crit}(f)$ are such that $f(x)>0$ and $f(y)<0$ and $i={\rm ind}(x)= {\rm
ind}(y)+1$.  Then
$$n^{f,\alpha}(x,y)~=~n^{F,\alpha}(x',y')~=~n^{G,\alpha'}(x',y')~.$$
All the counting of flow lines below will be done with respect to the
metric $\alpha'$ which we shall omit from the terminology.

As $\overline{G}$ is a chain map and in view of
the definition of the Morse complex we have
$$n^G(x',y')n^G(y',y)+\sum\limits_{z\in {\rm Crit}_i(G), z\not= y'}
n^G(x',z)n^G(z,y)~=~0~.$$
Now Lemma \ref{setting} implies that:
\begin{itemize}
\item[] if $z\in M\times \{1\}$, $z\not= y'$ then $n^G(z,y)=0$;
\item[] if $z\in (M\times \{0\})\backslash (N_{1/(2n)}\bigcup
N_{-1/(2n)})$ then $n^G(x',z)=0$ except for
$z=x$ and then $n^G(x',x)=1$;
\item[] if $z\in
N_{-1/(2n)}\bigcup N_{1/(2n)}$ then
$n^G(y',z)=0$;
\item[] if $z\in N_{1/(2n)}$ then $n^G(x',z)=0$.
\end{itemize}
As also $n^G(y',y)=1$, $n^{G}(x,y)=0$ we obtain
$$n^G(x',y')+\sum\limits_{z\in N_{-1/(2n)}} n^G(x',z)n^G(z,y)~=~0~.$$
Moreover for each $z\in N_{-1/2n}$, we also have
$$n^G(x',z)n^G(z,\theta z)+\sum\limits_{v\in {\rm Crit}_i(G)}
n^G(x',v)n^G(v,\theta z)~=~0$$ which gives
$$n^G(x',z)+n^G(x,\theta z)~=~0~.$$
So we get
$$\begin{array}{ll}
n^f(x,y)&=~n^G(x',y')\\[1ex]
&=~\sum\limits_{\overline{z}\in {\rm
{\rm Crit}_{i-1}}(g)}n^G(z,y)n^G(x,\theta z)\\[1ex]
&=~\sum\limits_{\overline{z}\in {\rm {\rm Crit}_{i-
1}}(g)}n^h(z,y)n^h(x, \theta z)~.
\end{array}$$
We can now define
$$h'~=~h|_{f^{-1}(-\infty,-1/(2n)]}~,~h''~=~h|_{f^{-
1}[1/(2n),\infty)}$$
and modulo the obvious identification of $N_{-1/(2n)}$ with
$N_{1/(2n)}$ we obtain the desired formula.

It remains now to prove the Lemma.

\noindent{\em Proof of Lemma \ref{deform}.}
We fix $p\in {\rm Crit}_i(F)\bigcap M\times \{1\}$.  We have
$p=p_1\times \{1\}$ with $p_1\in {\rm Crit}_{i-1}(F)\bigcap
M\times\{0\}$.
Consider $q\in {\rm Crit}_{i-1}(F)\bigcap M\times \{0\}$.
We have $n^F(p,q)=0$ if $q\not= p_1$ and $n^F(p,q)=1$ if $q=p_1$.
We shall now consider two cases:

(a)  $F(p)<0$.  In this case, by choosing $n$ sufficiently large, we
have
$$F(p)~=~F^{1/n}(p)\leq {\rm inf} (F^{1/n}(M^{1/n}))$$
and therefore $W^u_{F^{1/n}}\bigcap M^{1/n}=\emptyset$, so that
\begin{itemize}
\item[] $n^{F^{1/n}}(p,q)=n^F(p,q)$,
\item[] $n^{F^{1/n}}(p,r)=0$ for all $r\in N_{-1/(2n)}\bigcup
N_{1/(2n)}$.
\end{itemize}
As the number of critical points of $f$ is finite this implies that for
a sufficiently large $n$ the three properties of Lemma \ref{setting}
are verified for all critical points $p$ such that $F(p)<0$.

(b)  $F(p)>0$.  The same argument as above
shows that for all $q$ with $f(q)>0$ we have
that for sufficiently large $n$
$$n^{F^{1/n}}(p,q)~=~n^F(p,q)~.$$
There are two other properties that we still have to verify (and it is
only here that we shall need to use the specific form of $F^{1/n}$)\,:
\begin{itemize}
\item[(1)] if $f(q)<0$, then $n^{F^{1/n}}(p,q)=0$,
\item[(2)] if $r\in N_{1/(2n)}$, then $n^{F^{1/n}}(p,r)=0$.
\end{itemize}

We first note that for each $\epsilon$ the boundary of $M^{\epsilon}$
has the property that $W^{(-\epsilon,\epsilon)}\times
\{5\epsilon^{1/3}\}\bigcup N_{\epsilon}\times [0,5\epsilon^{1/3}]$ is
an entrance set into $M^{\epsilon}$ for the flow induced by the
negative of the gradient of $F^{\epsilon}$.  Similarly,
$N_{-\epsilon}\times [0,5\epsilon^{1/3}]$ is an exit set.  Assume that
for arbitrarily large $n$ the property (1) above is not satisfied.  As
the number of critical points of $f$ is finite this means that there is
a critical point $q$ of $f$ of index $i-1$ and a sequence
$n_k\rightarrow \infty$ such that $n^{F^{1/n_k}}(p,q)\not=0$ with
$f(q)<0$.  Clearly, as $n^F(p,q)=0$ this means that for each $n_k$
there is at least one flow line $\lambda_k$ (for the flow induced by
$F^{1/n_k}$) that joins $p$ to $q$ and that intersects $M^{1/n_k}$.
Let $x_k$ be the entrance point of this flow line into $M^{1/n_k}$ and
let $y_k$ be the exit point of this same flow line.  Inside $M^{1/n_k}$
the flow line joining $x_k$ to $y_k$ projects onto $N=N_0$ as a flow
line $\lambda'_k$ of $g$ joining the projection of $x_k$ denoted by
$x'_k$ to the projection $y'_k$ of $y_k$.  This happens because the
gradient of $F^{1/n_k}$ can be decomposed into two components, one
orthogonal to $N$ and the other having the same direction as the
gradient of $g$.  We now make $k\to \infty$.  Then, by compactness, we
may assume that $x_k\to x_0\in N$, $y_k\to y_0\in N$.  As in the
exterior of $M^{1/n_k}$ the function $F^{1/n_k}$ coincides with $F$ we
obtain that $y_0\in\overline{W}^s_F(q)$ and
$x_0\in\overline{W}^u_F(p)$.
 Because $F$ is a linear cobordism this
implies $x_0\in\overline{W}^u_f(p_1)$ and $y_0\in\overline{W}^s_f(q)$
(see Remark \ref{rem:class} (b)).
Of course, we also have $x'_k\to x_0$ and $y'_k\to y_0$.  Property (i)
of the function $g$ implies that $x_0\in \bigcup_{z\in {\rm Crit}_{\leq
i-2}} W^u_g(z)$.  As $x'_k$ and $y'_k$ were joined by the flow lines
$\lambda'_k$ of the flow induced by the negative of the gradient of $g$
we obtain that $y_0$ also belongs to $\bigcup_{z\in {\rm Crit}_{\leq
i-2}} W^u_g(z)$.  But this means that the intersection of some
$g$-unstable manifold of dimension $\leq i-2$ with the $f$-stable
manifold of $q$ is non trivial.  By property (ii) of $g$ that means
$${\rm dim}(M)-{\rm ind}_f(q)+(i-2)~=~{\rm dim}(W^s_f(q))+(i-1)\geq
{\rm dim}(M)~.$$
This implies ${\rm ind}(q)\leq i-2$ and leads to a contradiction.

We use a similar argument to deal with condition (2).  As above, if for
arbitrarily large $n$ (2) is not satisfied there is a critical point of
$F^{1/n_k}$, $r_k\in N_{1/(2n_k)}$ such that there exists at least one
$F^{1/n_k}$-flow line $\lambda_k$ joining $p$ to $r_k$ and ${\rm
ind}(r_k)=i-1$.  Obviously, when $k\to \infty$ we may assume that
$r_k\to r$ with $r$ a critical point of index $i-1$ of $g$.  As before,
let $x_k$ be the entrance point of $\lambda_k$ inside $M^{1/n_k}$ and
let $x'_k$ be the projection of $x_k$ onto $N$.  The portion of
$\lambda_k$ inside of $M^{1/n_k}$ projects to a $g$-flow line
$\lambda'_k$ joining $x'_k$ to $r$.  Let
$$x_0~=~{\rm lim}(x_k)~=~{\rm lim}(x'_k)~.$$
As above we obtain
$$x_0 \in\overline{W}^u_F(p)\bigcap
M\times\{0\}=\overline{W}^u_f(p_1)~.$$
At the same time we also have $x_0\in \overline{W^s_g(r)}$.
Property (ii) of $g$ implies
$${\rm dim}(W^s_g(r))+{\rm dim}(W^u_f(p_1))\geq n$$
which means
$$n-1-{\rm ind}_g(r)+{\rm ind}_f(p_1)\geq n~.$$
It follows that ${\rm ind}_g(r)\leq i-2$, which contradicts the
hypothesis on $r$ and completes the proof of the Lemma (we have
actually proved more: not only the relevant numbers of flow lines
counted with signs are zero but actually the respective sets of flow
lines are empty), and also of the gluing Theorem \ref{theo:glueadapt}.
\end{proof}

\section{Applications}\label{applications}

\subsection{Rigidity of the Floer complex}\label{Floer_rig}

The purpose of this section is to adapt the technique described in
\S \ref{subsubsec:rigid} to show a rigidity result for the Floer
complex
which is similar to Theorem \ref{theo:rig}. This adaptation is rather
immediate
as it only makes use of very standard
tools from the construction of Floer homology.

\subsubsection{The Floer complex.}\label{subsubsec:rec_floer}
We start by recalling - sketchily - some elements of the construction
of the Floer complex.
This construction is presented in detail in many sources and we
refer to them for details  \cite{Sal}, \cite{HoZe} \cite{SaZe}. We
follow
closely the last two sources mentioned as our sign conventions are like
there.

Let $M$ be a symplectic closed manifold with symplectic form
$\omega$ and of dimension $2m$. We shall work under the assumption
that $\omega$ vanishes when evaluated on $\pi_{2}(M)$. There exist
almost complex structures $J$ on $M$ that are compatible with
$\omega$ in the sense that the bilinear form $g_{J}(X,Y)=\omega(X,
J(x)Y)$, $X,Y\in T_{x}M$, gives a Riemannian metric on $M$. In
fact, the set of such compatible almost complex structures - which
will be denoted by $\mathcal{J}$ -  is contractible and therefore
there exists a well defined Chern class $c_{1}\in H^{2}(M;\ZZ)$.
We shall also assume $c_{1}|_{\pi_{2}(M)}=0$.

Fix a $1$-periodic hamiltonian $H$ on $M$. In other words, this is
a smooth function $H:S^{1}\times M\to \RR$ where $S^{1}=\RR/\ZZ$.
Such a hamiltonian induces a one parameter family of Hamiltonian
vector fields $X_{H}^{t}$, $t\in S^{1}$, on $M$ which is defined
by the equality $\omega(X_{H}^{t},Y)=-d(H_{t})(Y)$ that is
required to hold for all vector fields $Y$ on $M$. To this family
we associate the time dependent differential equation
\begin{equation}\label{eq:period}
\dot{x}(t)=X^{t}_{H}(x(t))~.~
\end{equation}
The solutions to this equation define a one-parameter family of
symplectic diffeomorphisms $\phi^{t}$
on $M$ given by $\phi^{t}_{H}(x(0))=x(t)$.

The problem which is addressed by the Floer complex is to estimate the
number of
solutions of equation (\ref{eq:period}) that satisfy $x(1)=x(0)$  and
are contractible
as loops in $M$. We denote by $\Lambda_{0}$ the space of all smooth,
contractible loops
in $M$ and we denote by $P(H)\subset \Lambda_{0}$ these $1$-periodic
solutions.

Let $a_{H}:\Lambda_{0}\to \RR$ be the function defined by
$$a_{H}(x)=-\int_{D}\overline{x}^{\ast}\omega +\int_{0}^{1}H(t,x(t))dt
~.~$$

Here $\overline{x}:D\to M$ is an extension of $x:S^{1}\to M$ to the
unit disk $D\subset \RR^{2}$,
$S^{1}=\partial D$. Because $\omega|_{\pi_{2}(M)}=0$, the first
integral in the formula
of $a_{H}$ is independent of the choice of the extension $\overline{x}$
and
$a_{H}$ is well defined.

The critical points of $a_{H}$ are identified with the elements of
$P(H)$ and
Floer's theory is a sort of Morse theory for functionals of this type.
Of course, a number of non-degeneracy conditions are needed.
The first - which we shall assume from now on - is that all elements of
$P(H)$
are required to be non-degenerate in the sense that
\begin{equation}\label{equ:nondeg1}
det(1-d\phi^{1}_{H}(x(0)))\not=0 ~.~
\end{equation}
This condition is generically satisfied and is
similar to the requirement that a function be Morse in the classical
situation. Given that
$M$ is compact and non-degenerate periodic orbits are isolated,
it follows that the number of elements of $P(H)$ is finite. To each
element
$x\in P(H)$ we may associate an integer $\mu_{H}(x)=\mu_{CZ}(x')$ which
will play the
role of the Morse index in the classical Morse case. Here
$\mu_{CZ}(x')$
is the Conley-Zehnder  index \cite{SaZe} of the path $x':[0,1]\to
Sp(2m)$ defined as follows.
We first choose $\overline{x}:D\to M$ extending $x$. We fix a
trivialization of the tangent bundle
of $M$ over $Im(\overline{x})$ and we use it to view the family
$d(\phi^{t}_{H})_{x(0)}$ as a
path $x':[0,1]\to Sp(2m)$ with $x'(0)=Id$ and $det(Id-x'(1))\not=0$
(because $x$ is non-degenerate).
One important point here is that because $c_{1}|_{\pi_{2}(M)}=0$ any
two
extensions $\overline{x}$ produce homotopic paths $x'$ and thus the
definition
of $\mu_{H}(x)$ is independent of the choice of extension.

The fundamental problem is that the gradient of the functional $a_{H}$
does not define a
flow. However, the solutions $u(s,t):\RR\times \RR/\ZZ\to M$ of the
equation
\begin{equation}\label{equ:orbits}
\frac{\partial u}{\partial s}+J(u)\frac{\partial u}{\partial t}+\nabla
_{x}H(t,u)=0
\end{equation}
may be viewed as ``negative gradient flow lines" of this functional.
Here $J\in \mathcal{J}$ and $\nabla_{x} H$ is the gradient of $H(t,-)$
with respect to the metric $g_{J}$.

The energy of a map $u:\RR\times S^{1}\to M$ is defined by
\begin{equation}\label{equ:energy}
E_{H}(u)=\frac{1}{2}\int_{-\infty}^{\infty}\int_{0}^{1}(|\frac{\partial
u}{\partial s}|^{2}+
|\frac{\partial u}{\partial t}-X_{H}^{t}(u)|^{2})dtds~.~
\end{equation}

We shall denote by $\mathcal{M}(H,J)$ the space of solutions $u$ of the
equation \ref{equ:orbits} that
verify $E_{H}(u)<\infty$. This space has the remarkable property to be
compact
and to decompose as the union of the spaces
\begin{eqnarray}\nonumber \mathcal{M}(y,x;H,J)=
 \{u\in \mathcal{M}(H,J) :
  \ \lim_{s\to-\infty}u(s,t)=y(t)\ , \ \lim_{s\to +\infty}u(s,t)=x(t)\}
\end{eqnarray}
where $x, y \in P(H)$.
Here, convergence is uniform in $t$ and, as it should happen for true
gradient flow lines of $a_{H}$,
the elements of $\mathcal{M}(y,x;H,J)$ also verify
$\lim_{s\to \pm\infty}\frac{\partial u}{\partial s}(s,t)=0$ uniformly
in $t$.
Moreover, for $u\in\mathcal{M}(y,x;H,J)$ we have $E_{H}(u)=a_{H}(y)-
a_{H}(x)$.
There is an obvious action of $\RR$ on the space $\mathcal{M}(y,x;H,J)$
and we shall denote
by $\mathcal{M}'(y,x;H,J)$ the resulting orbit space.

For our fixed almost complex structure $J$ there exists a subset
$\mathcal{H}_{reg}(J)$ of the second
Baire category inside $C^{\infty}(S^{1}\times M)$
such that whenever $H\in\mathcal{H}_{reg}(J)$
the natural linearization $D_{u,H,J}$ of the operator
$\overline{\partial}_{H,J}(u)=\frac{\partial u}{\partial s}+
J(u)\frac{\partial u}{\partial t}+\nabla _{x}H(t,u)$
is surjective at each $u$ with $\overline{\partial}_{H,J}(u)=0$
\cite{FlHoSa}.
We shall call the pair $(H,J)$ regular whenever $H$ satisfies
(\ref{equ:nondeg1}) and $H\in \mathcal{H}_{reg}(J)$. In the Morse case,
the analogue
of the condition $H\in \mathcal{H}_{reg}(J)$ is the Morse-Smale
transversality requirement (when
the Riemannian metric is fixed).

If $(H,J)$ is a regular pair, the surjectivity of $D_{u,H,J}$ for all
solutions $u$ of (\ref{equ:orbits})
implies that $\mathcal{M}(y,x;H,J)$ are manifolds whose dimension can
be computed to be
$\mu_{H}(y)-\mu_{H}(x)$. Therefore, the spaces $\mathcal{M}'(y,x;H,J)$
are manifolds of dimension
$\mu_{H}(y)-\mu_{H}(x)-1$. Moreover, in our setting, these manifolds
admit coherent orientations \cite{FlHo1}.
The spaces $\mathcal{M}'(y,x;H,J)$ also have natural compactifications
$\overline{\mathcal{M}'}(y,x;H,J)$
which are manifolds with boundary and
whose combinatorics is quite similar to that of the moduli spaces of
flow lines in the Morse-Smale case.
In particular, if $\mu_{H}(x)=\mu_{H}(y)-1$, then
$\mathcal{M}'(y,x;H,J)$ is finite and if
$\mu_{H}(x)=\mu_{H}(y)-2$, then we have the (oriented) equality:
$$\partial\overline{\mathcal{M}'}(y,x;H,J)=\coprod_{\stackrel{z\in
P(H)}{\mu_{H}(z)=\mu_{H}(y)-1}}
\mathcal{M}'(y,z;H,J)\times \mathcal{M}'(z,x;H,J)~.~$$

The Floer complex of $(H,J)$ is now defined by:
\begin{eqnarray}\nonumber CF_{i}(H,J)=\ZZ[x\in P(H) : \mu_{H}(x)=i],\ \
\nonumber d:CF_{i}(H,J)\to CF_{i-1}(H,J), \\
\nonumber d(x)=\sum_{z\in P(H),\mu_{H}(z)=\mu_{H}-1}n^{H,J}(x,z)
\end{eqnarray}
where $n^{H,J}(x,z)$ is the number of elements (counted with signs) of
$\mathcal{M}'(x,z;H,J)$.

\subsubsection{Rigidity.}\label{subsubsec:fl_rig} For a fixed
hamiltonian $H:S^{1}\times M\to \RR$
let $$\delta_{H}=\inf\{a_{H}(y)-a_{H}(x) : \ x,y\in P(H)\ ,\
\mathcal{M}(y,x;H,J)\not=\emptyset\}~.~$$
When all the $1$-periodic orbits of $X_{H}$ are non-degenerate the set
$P(H)$ is finite
and we have $\delta>0$.
The purpose of this section is to prove the following result.

\begin{theo}\label{theo:floer} Fix a regular pair $(H,J)$. If $(H',J')$
is another regular pair
with $||H'-H||_{0}\leq \delta_{H}/4$, then
the Floer complex $CF_{\ast}(H,J)$ is a retract of the complex
$CF_{\ast}(H',J')$.
\end{theo}

As in the Morse case, we immediately deduce.

\begin{cor}\label{cor:floer} Given a Hamiltonian $H:M\to \RR$ let
$J,J'$ be two
almost complex structures such that the pairs $(H,J)$
and $(H,J')$ are regular. There exists an isomorphism of chain
complexes
 $CF_{\ast}(H,J)\approx CF_{\ast}(H,J')$.
\end{cor}

\begin{proof} The proof of the theorem is based on the idea described
in \ref{subsubsec:rigid}
and, besides the usual approach to comparing the Floer complexes of two
different
regular pairs, it makes use of {\em monotone} homotopies as introduced
in \cite{FlHo}, \cite{CiFlHo}.

The standard way to compare the Floer complexes associated
to the two regular pairs $(H^{0},J^{0})$ and $(H^{1},J^{1})$ is as
follows
\cite{SaZe}\cite{HoZe}. Take smooth homotopies $H^{01}:\RR\times
S^{1}\times M\to \RR$
and $J^{01}: \RR\times M\to End(TM)$,  $J^{01}_{s}\in\mathcal{J},
\forall s\in\RR$ such
that there exists $R>0$ with the property that, for $s\geq R$, we
have $(H^{01}_{s}(x),J^{01}_{s}(x))=(H^{0}(x),J^{0}(x))$
and for $s\leq -R$, $(H^{01}_{s},J^{01}_{s})=(H^{1}(x),J^{1}(x))$,
$\forall x\in M$. Consider the equation:
\begin{equation}\label{equ:comp_fl}
\frac{\partial u}{\partial s}+J^{01}(s,u)\frac{\partial u}{\partial
t}+\nabla^{s}_{x}H^{01}(s,t,u)=0
\end{equation}
where $\nabla^{s}_{x}H^{01}(s,t,-)$ is the gradient of the
function $H(s,t,-)$ with respect to the Riemannian metric induced
by $J^{01}_{s}$.  Let $\mathcal{M}(H^{01},J^{01})$ be the space of
all finite energy solutions $u:\RR\times \RR/\ZZ\to M$ of this
equation (where the energy $E_{H^{01}}$ is defined by using
$H^{01}$ in (\ref{equ:energy}) instead of $H$). Again this space
is compact and it decomposes as the union of the spaces
$\mathcal{M}(y,x;H^{01},J^{01})$ containing those solutions $u$
that also satisfy
$$\lim_{s\to-\infty}u(s,t)=y(t),\ \lim_{s\to+\infty}u(s,t)=x(t)$$ for
some
with $x\in P(H^{0})$ and $y\in P(H^{1})$. Moreover, for generic choices
of $(H^{01},J^{01})$,
the associated linearized operator $D_{u,H^{01},J^{01}}$ is surjective
for each solution $u$ of
(\ref{equ:comp_fl}). Whenever $(H^{i},J^{i})$ are regular pairs for
$i\in\{0,1\}$ and this additional
condition is satisfied we shall say that the pair $(H^{01},J^{01})$
is a regular homotopy. As in the case of equation (\ref{equ:orbits}) a
consequence of regularity is that the spaces
$\mathcal{M}(y,x;H^{01},J^{01})$ are manifolds and they also admit
coherent orientations.
The dimension of such a manifold is given by $\mu_{H^{1}}(y)-
\mu_{H^{0}}(x)$. Define now
$\psi^{H^{01},J^{01}}:CF_{\ast}(H^{1},J^{1})\to CF_{\ast}(H^{0},J^{0})$
by
$$\psi^{H^{01},J^{01}}(y)=\sum_{x\in P(H^{0}),
\mu_{H^{0}}(x)=\mu_{H^{1}}(y)}n^{H^{01},J^{01}}(y,x)$$
where $n^{H^{01},J^{01}}(y,x)$ is the number of elements  in
$\mathcal{M}(y,x;H^{01},J^{01})$ (counted with signs). Obviously, the
key point is that
this application is a chain homomorphism which induces an isomorphism
in homology (\cite{SaZe},
\cite{Sal}).
For such a homotopy $H^{01}$ let
$$a_{H^{01}}(s,x)=-\int_{D}\overline{x}^{\ast}\omega
+\int_{0}^{1}H^{01}(s,t,x(t))dt ~.~$$
We now consider $u:\RR\to \Lambda_{0}$, $u\in
\mathcal{M}(y,x,H^{01},J^{01})$ and by composition
define $a(s)=a_{H^{01}}(s,u(s))$.
We derive $a(s)$ with respect to $s$ and we obtain (as, for example, in
\cite{FlHo}
p. 50 - except for a change in signs)
$$\frac{da}{ds}(s)= d(a_{H^{01}_{s}})(u(s,t))\frac{\partial u}{\partial
s}(s,t) + \int_{0}^{1}
\frac{\partial H^{01}}{\partial s}(s,t,u(s,t))dt$$ where
$d(a_{H})(\gamma)\xi$
is the derivative of the functional $a_{H}:\Lambda_{0}\to \RR$ at
$\gamma\in \Lambda_{0}$
in the direction of the vector field $\xi$ defined along $\gamma$. Now
\begin{eqnarray}\label{equ:derv_action}
d(a_{H^{01}_{s}})(u(s,t))\frac{\partial u}{\partial s}(s,t)=\\
\nonumber =\int_{0}^{1} g_{J^{01}_{s}}(J^{01}(s,u(s,t))\frac{\partial
u}{\partial t}(s,t) +\nabla^{s}_{x}H^{01}(s,t,u(s,t)),
\frac{\partial u}{\partial s}(s,t))dt
\end{eqnarray}
and, because $u$ verifies equation (\ref{equ:comp_fl}), this is always
negative.

If it happens that $H^{1}(t,x)>H^{0}(t,x)$ for all $(t,x)\in
S^{1}\times M$
we may take $H^{10}$ to be a monotone homotopy in the sense
that $\frac{\partial H^{01}}{\partial s}(s,t,x)\leq 0$ for all
$(s,t,x)\in \RR\times S^{1}\times M$
(this is quite analogue to our simple morse cobordisms from Definition
\ref{simple}).
These homotopies have been introduced and used by Floer and Hofer in
\cite{FlHo}.
If $H^{01}$ is monotone, the formula above shows
that $\frac{da}{ds}\leq 0$ which  shows that
\begin{equation}\label{equ:action_ineq}a_{H^{1}}(y)=a_{H^{10}}(-
R,y)\geq
a_{H^{10}}(R, x)=a_{H^{0}}(x)~.~
\end{equation}

Suppose now that we have two regular homotopies
$(H^{i(i+1)},J^{i(i+1)})$, $i\in\{0,1\}$
relating the regular pairs $(H^{i},J^{i})$ and $(H^{i+1},J^{i+1})$. It
is then possible to
glue the two homotopies together thus getting a new homotopy
$(H^{02},J^{02})$ which is
defined  by $H^{02}_{T}(s,t,x)=H^{01}(s+T,t,x)$ for $s\leq 0$,
$H^{02}_{T}(s,t,x)=H^{12}(s-T,t,x)$ for $s\geq 0$ and by analogue
formulae for $J^{02}$.
For $T$ sufficiently big this homotopy is regular and it verifies
$\psi^{H^{02},J^{02}}=\psi^{H^{01},J^{01}}\circ\psi^{H^{12},J^{12}}$
(\cite{SaZe}, Lemma 6.4; \cite{Sal} Lemma 3.11). Notice also that if
both $H^{01}$ and
$H^{12}$ are monotone, then so is $H^{02}$.

We now turn to the statement of the theorem. We fix the regular pair
$(H,J)$ and we let
$(H',J')$ be a regular pair with $|H'(t,x)-H(t,x)|\leq \delta_{H}/4$,
$\forall (t,x)\in S^{1}\times M$.
We have $H(t,x)-\delta_{H}/3< H'(t,x)< H(t,x)+\delta_{H}/3$, $\forall
(t,x)\in
S^{1}\times M$. Let $(H^{2},J^{2})=(H+\delta_{H}/3,J)$,
$(H^{1},J^{1})=(H',J')$,
$(H^{0},J^{0})=(H-\delta_{H}/3,J)$. We may find
monotone homotopies $(H^{01},J^{01})$, $(H^{12},J^{12})$ relating these
pairs such that
for all $s\in\RR$ we have $H^{0}\leq H^{01}_{s}\leq H^{1}\leq
H^{12}_{s}\leq H^{2}$.
For fixed $H^{ij}$ the pair $(H^{ij},J^{ij})$ is
regular for a generic choice of $J^{ij}$ inside the set of time-
dependent, period
one, almost complex structures compatible with $\omega$, \cite{FlHo}
Theorem 23 (the use
of time dependent almost complex structures does not change the form or
the behaviour
of any of the previous equations).

We consider
the ``glued" regular homotopy $(H^{02},J^{02})=(H^{02}_{T},J^{02}_{T})$
defined as above which is
monotone and regular for large $T$.

In essence, we have now constructed the two non-trivial faces of the
cube in the proof of
Theorem \ref{theo:rig} and, to end the proof, we need to show that
$\psi^{H^{02},J^{02}}$
is an isomorphism. We consider the monotone homotopy
$(\overline{H}^{02},J)$
such that $H^{02}(s,t,x)=H(t,x)+h(s)$ with $h:\RR\to \RR$ some
convenient, decreasing function
and $J_{s}=J$. Obviously, this homotopy is regular and
$\psi^{\overline{H}^{02},J}=id$.
Let $(G,\overline{J})$ be a smooth homotopy of homotopies
$G:[0,1]\times \RR\times S^{1}\times M\to \RR$
such that $G_{\lambda}=G(\lambda,-,-,-)$ verifies $G_{0}=H^{02}$,
$G_{1}=\overline{H}^{02}$
and we have similar identities for $\overline{J}$.  There exists a
notion
of regularity for such a homotopy of homotopies $(G,\overline{J})$. If
our
$(G,\overline{J})$ is regular in this sense it induces a chain homotopy
$\Psi=\psi^{G,\overline{J}}:CF_{\ast}(H^{2},J^{2})\to
CF_{\ast+1}(H^{0},J^{0})$ defined by
$$\psi^{G,\overline{J}}(y)=\sum_{x\in
P(H^{0}),\mu_{H^{0}}(x)=\mu_{H^{2}}(y)+1}n^{G,\overline{J}}(y,x)$$
where $n^{G,\overline{J}}(y,x)$ is the number of elements (counted with
signs) of the
moduli spaces
$$\mathcal{M}(y,x;G,\overline{J})=\{(u,\lambda) : \ \lambda\in [0,1],\
u\in\mathcal{M}(y,x;
G_{\lambda},\overline{J}_{\lambda}) \}$$
which, are finite by regularity  when
$\mu_{H^{0}}(x)=\mu_{H^{2}}(y)+1$ and also admit a coherent choice of
orientations
(see \cite{SaZe} Lemma 6.3, \cite{FlHo} \S 4.3). With this definition
we have $d\Psi-\Psi d= \psi^{\overline{H}^{02},J}-
\psi^{H^{02},J^{02}}=Id-\psi^{H^{02},J^{02}}$.
We index the elements in $P(H^{2})=P(H^{0})=P(H)=\{p_{1},p_{2},\ldots
p_{k}\}$
such that $a_{H}(p_{i})< a_{H}(p_{j}) \Rightarrow i < j$
and we denote by $A=(a_{ij})$ the matrix of the linear map $d\Psi -\Psi
d$ in this basis.

Because both $H^{02}$ and $\overline{H}^{02}$ are monotone homotopies
we may take
$G$ to be a homotopy of monotone homotopies in the sense that
$G_{\lambda}$ is a monotone
homotopy for all $\lambda\in [0,1]$. Simultaneously, we also need the
pair $(G,\overline{J})$
to be regular. This can be achieved generically (with again
$\overline{J}$ time-dependent)
as discussed in \cite{FlHo} \S 4.3. In this case, it follows from
(\ref{equ:action_ineq}) that
$n^{G,\overline{J}}(y,x)\not=0 \ \Rightarrow a_{H^{2}}(y)\geq
a_{H^{0}}(x)$. The coefficients of the
Floer differential in $CF_{\ast}(H,J)$ verify $n^{H,J}(y,x)\not=0 \
\Rightarrow
a_{H}(y)\geq \delta_{H}+a_{H}(x)$. By looking to the matrix $A$ we
obtain
that $a_{ij}\not=0 \ \Rightarrow a_{H^{2}}(p_{j}) \geq
a_{H^{0}}(p_{i})+\delta_{H}$.
But $a_{H}(p_{j})+\delta_{H}/3= a_{H^{2}}(p_{j})$,
$a_{H^{0}}(p_{i})=a_{H}(p_{i})-\delta_{H}/3$
so we get that if $a_{ij}\not=0$, then $a_{H}(p_{j})> a_{H}(p_{i})$
which means $j>i$. Therefore,
$A$ is upper triangular and we conclude that $\psi^{H^{02},J^{02}}$ is
an isomorphism which ends the proof.
\end{proof}

\begin{rem}\label{rem:ext}
{\rm It is clear that the result above can be extended in many ways.
One such possibility is to consider
monotone symplectic manifolds which means that the two morphisms
$c_{1}:\pi_{2}(M)\to \ZZ$
and $[\omega]:\pi_{2}(M)\to \RR$ are proportional. Besides other
technical complications, this implies that
the action-functional $a_{H}$ is $S^{1}$ valued. In the next section we
discuss rigidity (and gluing)
for $S^{1}$ valued Morse functions. It is quite likely that
the rigidity statement proved below in the $S^{1}$-valued case carries
over to the Floer monotone setting.
}\end{rem}

\subsection{Circle-valued functions}\label{Novikov}
The purpose of this subsection is to apply the methods of Section
\ref{Morse cobordism} to
$S^{1}$-valued functions. We shall obtain a rigidity statement similar
to that in
\ref{subsubsec:rigid} and we shall also adapt to this context the
gluing result
\ref{theo:glueadapt}. This is then used to show how to algebraically
recompose the Morse-Novikov
complex out of information on a single fundamental domain.

\subsubsection{The Novikov complex}\label{subsec:intro_nov}

We now consider a circle-valued Morse-Smale
function on a compact manifold $M$
$$(f,\alpha)~:~M \to S^1~=~\RR/\ZZ~.$$
In this section we fix the main definitions and properties
of the Novikov complex of $f$.

The pullback infinite cyclic cover
$$\overline{M}~=~f^*\RR~=~
\{(x,t) \in M \times \RR\,\vert\,f(x)=[t] \in
S^1\}$$ is non-compact, and $(f,\alpha)$ lifts
to a $\ZZ$-equivariant Morse-Smale function
$$(\overline{f},\alpha)~:~\overline{M}~=~f^*\RR \to \RR$$
with
$$\overline{f}~:~\overline{M}\to \RR~;~(x,t) \mapsto t~.$$
$$\xymatrix{\overline{M} \ar[r]^{\overline{f}} \ar[d] & \RR \ar[d] \\
M \ar[r]^f & S^1}$$ and where we identify the
metric $\alpha$ with its pull-back to
$\overline{M}$. The generating covering
translation along the downward gradient flow of
$\overline{f}$ is given by
$$z~:~\overline{M}\to \overline{M}~;~(x,t) \mapsto (x,t-1)~.$$
We shall assume that $M$ and $\overline{M}$ are
connected, so that there is defined an exact
sequence of groups
$$\xymatrix{\{1\} \ar[r] & \pi_1(\overline{M}) \ar[r] &
\pi_1(M) \ar[r]^-{\displaystyle{f_*}}
&\pi_1(S^1)=\ZZ \ar[r] & \{1\}}~.$$ Let
$\widetilde{M}$ be a regular cover of
$\overline{M}$ with group of covering
translations $\pi$ (e.g. $\overline{M}$ itself
with $\pi=\{1\}$, or the universal cover of $M$
with $\pi=\pi_1(\overline{M})$). Then
$\widetilde{M}$ is a regular cover of $M$, and
the group of covering translations $\Pi$ fits
into an exact sequence
$$\xymatrix{\{1\} \ar[r] & \pi \ar[r] &
\Pi \ar[r] &\ZZ \ar[r] & \{1\}~,}$$ so that
$$\Pi~=~\pi \times_{\zeta} \ZZ$$
with $\zeta$ the monodromy automorphism defined
by
$$\zeta~:~ \pi \to \pi~;~g \mapsto z^{-1} g z$$
for any lift of a generator $1 \in \ZZ$ to an
element $z \in \Pi$. The group ring
$$\ZZ[\Pi]~=~\ZZ[\pi]_{\zeta}[z,z^{-1}]$$
consists of the polynomials
$$\sum\limits^{\infty}_{j=-\infty}a_jz^j~~(a_j \in \ZZ[\pi])$$
with $\{j \in \ZZ\,:\,a_j \neq 0\}$ finite and
$$az~=~z\zeta(a)~~(a \in \ZZ[\pi])~.$$
For any $\ZZ[\pi]$-module $F$ and $j \in \ZZ$
let $z^jF$ be the $\ZZ[\pi]$-module with
elements $z^jx$ ($x \in F$) and
$$z^jx+z^jy~=~z^j(x+y)~,~az^jx~=~z^j\zeta^j(a)x \in z^jF~.$$
\indent Let
$$(\overline{f},\alpha)~:~\overline{M} \to \RR$$
be a lift of $(f,\alpha)$ to a
$\ZZ$-equivariant real-valued Morse function on
the infinite cyclic cover $\overline{M}$ of
$M$. Lift each critical point of $f:M \to S^1$
to a critical point of
$\overline{f}:\overline{M} \to \RR$, allowing
the identification
$${\rm Crit}_i(\overline{f})~=~\bigcup
\limits^{\infty}_{j=-\infty}z^j{\rm Crit}_i(f)~.$$

\begin{defi} {\rm
(i) The {\it Novikov ring}
$$\ZZ[\pi]_\zeta((z))~=~\ZZ[\pi]_{\zeta}[[z]][z^{-1}]$$
is the ring of formal power series
$$\sum\limits^{\infty}_{j=-\infty}a_jz^j~~(a_j \in \ZZ[\pi])$$
such that $\{j \leq 0\,:\,a_j \neq 0\}$ is finite.\\
(ii)  The {\it Novikov complex} of $(f,\alpha)$ is the based f.g. free
$\ZZ[\pi]_\zeta((z))$-module chain complex $C=C^{Nov}(M,f,\alpha)$ with
$$\begin{array}{l}
d~:~C_i~=~\ZZ[\pi]_\zeta((z))[{\rm Crit}_i(f)]\to
C_{i-1}~=~\ZZ[\pi]_\zeta((z))[{\rm Crit}_{i-1}(f)]~;\\[2ex]
\hspace*{20mm} x \mapsto
\sum\limits^{\infty}_{j=-\infty}
\sum\limits_{g\in \pi} \sum\limits_{y\in {\rm
Crit}_{i-1}(f)}n^{\widetilde{f},\alpha}(x,z^jgy)z^jgy~.
\end{array}$$
\hfill\qed}
\end{defi}

Assume that $1 \in S^1$ is a regular value of $f$ and let $N=f^{-1}(1)
\subset M$.  Cutting $M$ along $N$ gives a fundamental domain
$(M_N;N,z^{-1}N)$ for $\overline{M}$ with a Morse-Smale function
$$(f_N,\alpha_N)~:~(M_N;N,z^{-1}N) \to ([0,1];\{0\},\{1\})$$
which is constant on $N$ and $z^{-1}N$.
$$\xymatrix@C+2pt{\ar@{.}[r]&\ar@{-}[rrrrrr]
\ar@{-}[dd]_{\displaystyle{zN}} &&
\ar@{-}[dd]_{\displaystyle{N}} &&
\ar@{-}[dd]_{\displaystyle{z^{-1}N}}
&& \ar@{-}[dd]^{\displaystyle{z^{-2}N}}\ar@{.}[r]&\\
&& zM_N && ~M_N~ && z^{-1}M_N && \\
\ar@{.}[r]&\ar@{-}[rrrrrr] &&&&&&\ar@{.}[r]&}
$$
Each critical point $x \in M$ of $f$ has a
unique lift to a critical point of $f_N$, which
is also denoted by $x \in M_N$, with
$$\begin{array}{l}
\overline{f}~=~\bigcup\limits^{\infty}_{j=-\infty}z^jf_N~:~
\overline{M}~=~\bigcup\limits^{\infty}_{j=-\infty}z^jM_N
\to
\RR~,\\[4ex]
{\rm Crit}_i(f_N)~=~{\rm Crit}_i(f)~,~ {\rm
Crit}_i(\overline{f})~=~\bigcup\limits^{\infty}_{j=-\infty}z^j
{\rm Crit}_i(f_N)~,\\[2ex]
n^{\overline{f},\alpha}(z^jx,z^ky)~=~
\begin{cases}
n^{f,\alpha}(x,y)&\hbox{\rm if $j=k$}\\
0&\hbox{\rm if $j>k$}
\end{cases}~~(x \in {\rm Crit}_i(f),y \in {\rm Crit}_{i-1}(f))~.
\end{array}$$
Use the lifts of the critical points of $f$ to
$M_N$ in the construction of the Novikov
complex $C^{Nov}(M,f,\alpha)$, so that the
coefficients of the differentials are in
$\ZZ[\pi]_{\zeta}[[z]] \subset
\ZZ[\pi]_{\zeta}((z))$.

\begin{lem} \label{lem:Nov}
The Novikov complex of $(f,\alpha):M \to S^1$ is
$$C^{Nov}(M,f,\alpha)~=~\ZZ[\pi]_{\zeta}((z))\otimes_{\ZZ[\pi]_{\zeta}[
[z]]}
\invliml
C(M_N(\ell),f_N(\ell),\alpha_N(\ell))$$ with
$C(M_N(\ell),f_N(\ell),\alpha_N(\ell))$ the
$\ZZ[\pi]$-coefficient Morse complex of the
Morse-Smale function
$$(f_N(\ell),\alpha_N(\ell))~=~(\overline{f},\alpha)\vert~:~
(M_N(\ell);z^{\ell}N,z^{-1}N) \to \RR$$ on
$$M_N(\ell)~=~\bigcup^{\ell}_{j=0}z^jM_N \subset \overline{M}~=~
\bigcup^{\infty}_{j=-\infty}z^jM_N~.$$ The
inverse limit is a based f.g. free
$\ZZ[\pi]_{\zeta}[[z]]$-module chain complex,
where the limit is with respect to the natural
projections
$$C(M_N(\ell+1),f_N(\ell+1),\alpha_N(\ell+1))
\to C(M_N(\ell),f_N(\ell),\alpha_N(\ell))~.$$
\end{lem}
\begin{proof} Let $F_i$ be the based f.g. free $\ZZ[\pi]$-module
generated
by the index $i$ critical points of $f$. For
any $\ell \geq 0$
$$C(M_N(\ell),f_N(\ell),\alpha_N(\ell))_i~=~
\sum\limits^{\ell}_{j=0}z^jF_i~,$$ and
$$\invliml C(M_N(\ell),f_N(\ell),\alpha_N(\ell))_i~=~
\invliml
\sum\limits^{\ell}_{j=0}z^jF_i~=~(F_i)_{\zeta}[[z]]$$
so that
$$\begin{array}{ll}
C^{Nov}(M,f,\alpha)_i&=~(F_i)_{\zeta}((z))\\[1ex]
&=~\ZZ[\pi]_{\zeta}((z))\otimes_{\ZZ[\pi]_{\zeta}[[z]]}
\invliml
C(M_N(\ell),f_N(\ell),\alpha_N(\ell))_i~.
\end{array}$$
\end{proof}

\subsubsection{Assembly and rigidity for the Novikov complex}

It is clear that the Novikov complex of $(f,\alpha):M \to S^1$ only
depends on the geometric behaviour of the lift $\overline{f}$ on a
single fundamental domain.  With the notations of the last section, one
expects therefore that there should exist an assembly process that
would consist in extracting from $f_N$ (or, equivalently, from
$f|_{M\backslash U(N)}$ with $U(N)$ a small tubular neighbourhood of
$N$) precisely the data necessary to algebraically piece together
$C^{Nov}(M,f,\alpha)$.  The work of Pajitnov \cite{P} describes one
such method.  We shall present below a different one which is based on
the techniques described in the first parts of the paper.  We
will compare in more detail our results with those of Pajitnov at the
end of the section.

The notion of a splitting of a real-valued
Morse function has an obvious analogue for circle-valued Morse
functions.
As in the real valued case (see just above Definition \ref{splitting}),
there is no restriction to assume that $\alpha$ equals the product
metric
$\beta+dt^2$ in a tubular neighbourhood of $N$ with $\beta=\alpha|_N$.

\begin{defi} \label{splitting2} {\rm A $(\epsilon,\tau)$-{\it
splitting} of a
Morse-Smale function $(f,\alpha):M \to S^1$
along a Morse-Smale function
$(g,\beta):N=f^{-1}(1)\to \RR$ with
$\beta=\alpha\vert_N$ is a Morse-Smale function
$(h,\alpha):M \to S^1$ such that:
\begin{itemize}
\item[(i)] $h\vert_N =e^{2\pi i\tau} g:N \to J$ for a convenient
interval $J\subset S^1$.
\item[(ii)] $(h,\alpha)=(f,\alpha)$ except in a small tubular
neighbourhood
$W=N \times [-\epsilon,\epsilon] \subset M$ of
$N=N \times \{0\} \subset M$. We let
$M'=M\backslash N$.
\item[(iii)] With respect to this parametrization $f(x,t)=e^{2\pi i
t}\in J'$ for $(x,t)\in
N\times (-\epsilon,\epsilon)$ with $J'\subset
S^1$ an interval containing $1$ and
$$\begin{array}{l}
{\rm Crit}_i(h)~=~{\rm Crit}_i(f)\cup ({\rm
Crit}_{i-1}(g)\times \{-\epsilon/2\}) \cup
({\rm Crit}_i(g)\times \{\epsilon/2\})~.
\end{array}$$
\item[(iv)] The restrictions of $(h,\alpha)$ to the submanifolds
$$\begin{array}{l}
W_{\epsilon}~=~N \times [-\epsilon/2,\epsilon/2] \subset W~,\\[1ex]
M'_{\epsilon}~=~M' \backslash \big(N \times
(-\epsilon/2,\epsilon/2)\big) \subset
M'~,\\[1ex]
\end{array}$$
are Morse-Smale cobordisms
$$\begin{array}{l}
(g_{\epsilon},\beta_{\epsilon})~:~
(W_{\epsilon};N\times \{\epsilon/2\},N\times \{-\epsilon/2\}) \to
J'',\\[1ex]
(h',\alpha'):~
(M'_{\epsilon};N\times \{\epsilon/2\},N\times
\{-\epsilon/2\}) \to J'''\subset S^1-\{1\}.
\end{array}$$
with $J''\cup J'''=S^1$.
\item[(v)] The cobordism $(g_{\epsilon},\beta_{\epsilon})$ is linear,
so
that
$$C(W_{\epsilon},g_{\epsilon},\beta_{\epsilon})~=~C(1:C(N,g,\beta) \to
C(N,g,\beta))~. \eqno{\qed}$$
\end{itemize} }
\end{defi}

The key point of the construction is that, by contrast to the
Novikov complex of $(f,\alpha)$, the Novikov complex of a
splitting $(h,\alpha)$ is very simple because the flow lines of
the corresponding lift $\overline{h}$ are never longer than one
fundamental domain. Indeed, the only way a flow line of $h$ can
cross $N$ is by joining inside $N\times [-\epsilon/2,\epsilon/2]$
a critical point contained in $N\times \{-\epsilon/2\}$ to one in
$N\times \{\epsilon/2\}$.  As $h$ restricts to a linear cobordism
on $N\times [-\epsilon/2,\epsilon/2]$ this leads to a simple
description for the Novikov complex of $h$.

For any such splitting we may identify the interval $J'''$ with
$[\epsilon/2,1-\epsilon/2]$ (such that, in $S^1$,
$1-\epsilon=-\epsilon)$ and if we denote $N'=N\times \{\epsilon/2\}$
and $N'' =N\times \{1-\epsilon/2\}$, then the cobordism $(h',\alpha')$
equals $(\overline{h},\alpha)|_{\overline{h}^{-1}
([\epsilon/2,1-\epsilon/2])}$.  Clearly, $(M_{\epsilon}';N',N'')$ is
diffeomorphic to $(M_N;N,z^{-1}N)$.  The cobordism
$$(h',\alpha'_N)~=~(h_N,\alpha_N)\vert~:~(M'_{\epsilon};N',N'') \to
\RR$$
has
$$\begin{array}{l}
{\rm Crit}_i(h')~=~({\rm Crit}_i(g)\times\{\epsilon/2\}) \cup {\rm
Crit}_i(f_N) \cup ({\rm Crit}_{i-1}(g)\times \{1-\epsilon/2\})~.
\end{array}$$
\indent The algebraic data contained in a splitting of
an $S^1$-valued map is brought together in the following structure.

\begin{defi} \label{geocob2}
{\rm The {\it algebraic cobordism} of a splitting $(h,\alpha):M \to
S^1$ of
$(f,\alpha):M \to S^1$ along $g$ is defined by
$$\Gamma_N(M,h,\alpha)~=~(F,D,\theta,\theta',\psi)$$
with $\ZZ[\pi]$-module chain complexes and chain maps
$$\begin{array}{l}
D~=~C(N,g,\beta)~,~F~=~C(M'_N\backslash
\partial M'_N,h'_N\vert,\alpha'_N\vert)~=~
C(M_N,f_N,\alpha_N)~,\\[2ex]
\theta~:~F_i \to D_{i-1}~;~x \mapsto
\sum\limits_{y \in {\rm Crit}_{i-1}(g)}
n^{h'_N,\alpha'_N}(x,(y\times\{\epsilon/2\}))y~,\\[2ex]
\theta'~:~D_i \to F_i~;~x \mapsto
\sum\limits_{y \in {\rm Crit}_i(f_N)}
n^{h'_N,\alpha'_N}((x\times \{1-\epsilon/2\}),y)y~,\\[2ex]
\psi~:~D_i \to zD_i~;~x \mapsto \sum\limits_{y
\in {\rm Crit}_i(g)}
n^{h'_N,\alpha'_N}((x\times\{1-\epsilon/2\}),(y\times
\{\epsilon/2\}))y~.
\end{array}$$
\hfill\qed}
\end{defi}

By analogy with Proposition \ref{theta1} :

\begin{prop} \label{theta2}
Let $(f,\alpha):M \to S^1$ be a Morse-Smale function, with $0 \in S^1$
a regular value.\\
{\rm (i)} For every Morse-Smale function $(g,\beta):N=f^{-1}(0) \to
\RR$
with $\beta=\alpha\vert_N$ and every $\delta >0$ there exists an
$(\epsilon,\tau)$-splitting of $(h,\alpha)$ along $(g,\beta)$
such that $\Vert h - f \Vert_0 \leq \delta$.\\
{\rm (ii)} The Novikov complex of a splitting $(h,\alpha)$ of
$(f,\alpha)$
along $(g,\beta)$ with algebraic cobordism
$\Gamma_N(M,h,\alpha)=(F,D,\theta,\theta',\psi)$ is given by
$$C^{Nov}(M,h,\alpha)~=~E$$
with $E$ the based f.g. free $\ZZ[\pi]_{\zeta}((z))$-module chain
complex defined by
$$\begin{array}{l}
d_E~=~\begin{pmatrix} -d_D & 0 & 0 \cr
1-z\psi & d_D & \theta \cr -z\theta' & 0 &
d_F \end{pmatrix}~:\\[4ex]
E_i~=~(D_{i-1}\oplus D_i \oplus F_i)_{\zeta}((z))
\to E_{i-1}~=~(D_{i-2}\oplus D_{i-1} \oplus F_{i-1})_{\zeta}((z))~.
\end{array}$$
\end{prop}
\begin{proof}
(i) As for \ref{theta2} (i).\\
(ii) The Morse complex of
$(\overline{h},\overline{\alpha}):\overline{M} \to \RR$ is
$$\begin{array}{l}
d_{C(\overline{M},\overline{h},\overline{\alpha})}~=~
\begin{pmatrix} -d_D & 0 & 0 \cr
1-z\psi & d_D & \theta \cr -z\theta' & 0 &
d_F \end{pmatrix}~:\\[4ex]
C(\overline{M},\overline{h},\overline{\alpha})_i~=~
(D_{i-1}\oplus D_i \oplus F_i)_{\zeta}[z,z^{-1}]\\[2ex]
\hskip100pt \to C(\overline{M},\overline{h},\overline{\alpha})_{i-1}~=~
(D_{i-2}\oplus D_{i-1} \oplus F_{i-1})_{\zeta}[z,z^{-1}]
\end{array}$$
and
$$C^{Nov}(M,h,\alpha)~=~E~=~\ZZ[\pi]_{\zeta}((z))
\otimes_{\ZZ[\pi]_{\zeta}[z,z^{-1}]}
C(\overline{M},\overline{h},\overline{\alpha})~.$$
\end{proof}

In the setting of Proposition \ref{theta2} define
a based f.g. free $\ZZ[\pi]_{\zeta}((z))$-module chain
complex $\widehat{F(h)}$ by
$$\begin{array}{l}
d_{\widehat{F(h)}}~=~d_F+z\theta'(1-z\psi)^{-1}\theta~=~
d_F+\sum\limits^{\infty}_{j=0}z^{j+1}\theta'\psi^j\theta~:\\[1ex]
\hskip100pt \widehat{F(h)}_i~=~(F_i)_{\zeta}((z)) \to
\widehat{F(h)}_{i-1}~=~(F_{i-1})_{\zeta}((z))
\end{array}$$
The role of this complex is to eliminate in a purely algebraic way the
contribution of $g$ from the Novikov complex of the splitting
$(h,\alpha)$.
Clearly, this suggests that $\widehat{F(h)}$ is a good approximation of
$C^{Nov}(M,f,\alpha)$ and we shall see further that this is indeed the
case.

It is convenient to fix the notations
$$\mathcal{R}~=~\ZZ[\pi]_{\zeta}((z))~~,~~
\mathcal{R}'~=~\ZZ[\pi]_{\zeta}[[z]]\subset \mathcal{R}~.$$
All complexes and isomorphisms in this result are understood to belong
to the category of $\mathcal{R}$-module chain complexes.

\begin{defi} {\rm An isomorphism of based f.g. free $\mathcal{R}$-
modules
$\Theta:F \to G$ is  {\it ${\mathcal R}'$-simple} if
$$\tau(\Theta) \in {\rm im}(K_1(\mathcal{R}') \to
K_1(\mathcal{R}))~.\eqno{\qed}$$}
\end{defi}

The result below is the $S^1$-analogue of the Rigidity Theorem
\ref{theo:rig} (iii) and the gluing Theorem \ref{theo:glueadapt}.

\begin{theo} \label{theo:Novikov} Fix the Morse-Smale function
$(f,\alpha):M \to S^1$.\\
{\rm (i)}There exists $\delta_{f}>0$ such that if $(f',\alpha')$ is a
Morse-Smale function
with $||f'-f||_{0}\leq \delta_{f}$, then $C^{Nov}(M,f,\alpha)$ is a
retract of $C^{Nov}(M,f',\alpha')$.
In particular, the $\mathcal{R}'$-simple isomorphism type of the
Novikov complex
of $(f,\alpha)$ is independent of the choice of $\alpha$.\\
{\rm (ii)} For any $(\epsilon,\tau)$-splitting $(h,\alpha):M \to S^1$
of $(f,\alpha)$ along a Morse-Smale function $(g,\beta):N \to \RR$ with
$\epsilon,\tau$ sufficiently small there is an $\mathcal{R}'$-simple
isomorphism
$$\Theta ~:~C^{Nov}(M,f,\alpha) \xrightarrow[]{\cong}
\widehat{F(h)}~.$$
{\rm (iii)} For any $n\in {\bf N}$ there exist Morse-Smale functions
$(g_n,\beta): N\to \RR$ and splittings $(h_n,\alpha):M\to S^1$ of
$(f,\alpha)$ along $(g_n,\beta)$ such that there is an
$\mathcal{R}'$-simple isomorphism
$$\Theta^n~:~C^{Nov}(M,f,\alpha) \xrightarrow[]{\cong}
\widehat{F(h_n)}$$
with $\Theta^n(x)-x\in z^n\mathcal{R}'[{\rm Crit}_{\ast}(f)]$ for each
basis element $x$.\hfill$\qed$
\end{theo}

\begin{rem} \label{psi2} {\rm The complex $\widehat{F(h)}$ behaves as
if it were the Novikov complex $C^{Nov}(M,\widehat{f},\alpha)$ of a
Morse-Smale function $(\widehat{f},\alpha):M \to S^1$ obtained from a
splitting $(h,\alpha)$ by cancelling the pairs of critical points added
to those of $f$ in the construction of $h$, with
$${\rm Crit}_*(\widehat{f})~=~{\rm Crit}_*(f)~.$$
Indeed, as in the real-valued case (see Remark \ref{rem:glue1})
such a function does exist (this follows from Hutchings' extension
in \cite{H} of the bifurcation analysis of Laudenbach \cite{Lau})
and one may construct it such that it is $C^{0}$-close to $f$.
Using point (i) of the theorem and the existence of such an
$\widehat{f}$ is sufficient to prove (ii) but, again as in the
real-valued case, it is not enough for (iii).  In the following we
shall prove the whole statement independently of bifurcation
considerations.  \hfill$\qed$}
\end{rem}

\begin{proof}
Order the critical points of $f_N$ following
their critical values. Explicitly, we fix a
total order on the finite set
$S={\rm Crit}_{\ast}(f_N)={\rm Crit}_{\ast}(f)$ such that
$x\leq y$ implies $f_N(x)\leq f_N(y)$ for
$x,y\in S$.\\
\indent
The proof of the theorem uses the adaptation to the $S^1$-valued case
of the techniques developed in the first three sections of this paper.
These methods will produce certain morphisms relating the various chain
complexes involved.  Showing that these morphisms are isomorphisms is
less immediate than in the $\RR$-valued case, but it will always follow
from the fact that our morphisms satisfy the assumptions of the rather
obvious lemma below.

\begin{lem}\label{lem:algNov} Consider the free $\mathcal{R}$-module
$\mathcal{R}\{S\}$ generated by a finite ordered set
$S$.  An $\mathcal{R}$-module endomorphism $\Theta :
\mathcal{R}\{S\}\to
\mathcal{R}\{S\}$ such that for each $x\in S$
$$\Theta(x)-x \in  z\mathcal{R'}\{S\} + \mathcal{R'}\{ y \in S: y < x
\}$$
is an $\mathcal{R}'$-simple automorphism.
\end{lem}
\begin{proof} Write $S=\{x_1,x_2,\ldots, x_k\}$, with $x_i<x_j$ for
$i<j$. The matrix of $\Theta$ has entries in $\mathcal{R}'$, of the
form
$$\Theta~=~\begin{pmatrix} 1+a_{11} & b_{12} & b_{13} & \dots & b_{1k}
\\
a_{21} & 1+a_{22} & b_{23} & \dots & b_{2k} \\
a_{31} & a_{32} & 1+a_{33} & \dots & b_{3k} \\
\vdots & \vdots & \vdots & \ddots & \vdots \\
a_{k1} & a_{k2} & a_{k3} & \dots & 1+a_{kk} \end{pmatrix}$$
with $a_{ij} \in z\mathcal{R'}=z\ZZ[\pi]_{\zeta}[[z]]$,
 $b_{ij} \in \mathcal{R}'$. The augmentation
 $$\epsilon~:~\mathcal{R}' \to \ZZ[\pi]~;~z \mapsto 0$$
 sends $\Theta$ to an upper triangular matrix of the form
$$\epsilon(\Theta)~=~\begin{pmatrix} 1 & \epsilon(b_{12}) &
\epsilon(b_{13})
& \dots & \epsilon(b_{1k}) \\
0 & 1 & \epsilon(b_{23}) & \dots & \epsilon(b_{2k}) \\
0 & 0 & 1 & \dots & \epsilon(b_{3k}) \\
\vdots & \vdots & \vdots & \ddots & \vdots \\
0 & 0 & 0 & \dots & 1 \end{pmatrix}$$
which is clearly invertible. The matrix defined by
$$\Psi~=~\Theta\epsilon(\Theta)^{-1} -I $$
has entries in $z\mathcal{R}'$. Thus $\Theta$ itself is invertible,
with
inverse
$$\Theta^{-1}~=~\epsilon(\Theta)^{-1}(1+\sum\limits^{\infty}_{j=1}(-
\Psi)^j)$$
defined over ${\mathcal R}'$.
\end{proof}

We now return to the proof of Theorem \ref{theo:Novikov}.  The first
remark is that the Definitions \ref{Morse-Smale}, \ref{cobordism} and
\ref{simple} apply ad literam to the case of $S^1$-valued maps because
the relevant conditions are local in nature.  Similarly, the statements
in Proposition \ref{chain complex} (i) and Proposition \ref{chaineq}
(iii) remain valid if we replace the relevant Morse complexes by the
corresponding Morse-Novikov complexes with one exception: the argument
given for the proof of \ref{chaineq} (iii) only shows that $f^M$ is a
chain map in the $S^1$-valued case.

More care is necessary for the constructions of Morse cobordisms.

\begin{lem}\label{constr4} {\rm (i)} Let $(f_0,\alpha_0),
(f_1,\alpha_1): M\to S^1$ be two Morse-Smale functions with $M$
closed.  Suppose that $f_0$ and $f_1$ are homotopic.  Then there exists
$c\in \RR$ and a simple Morse cobordism $(F,\alpha'):M\times [0,1]\to
S^1$ between $(e^{2\pi i c}f_0,\alpha_0)$ and $(f_1,\alpha_1)$.\\
{\rm (ii)} If there exists a constant $c\geq 0$ such that
$$(f_1,\alpha_1)~=~(e^{2\pi i c}f_0,\alpha_0)~:~M \to S^1$$
then there exists a cobordism $(F,\alpha'):M\times [0,1]\to S^1$
between $f_0$ and $f_1$ such that $F^{M\times [0,1]}$ is the identity.
\end{lem}
\begin{proof} The difference between proving this statement and the
constructions in \S\ref{subsubsec:constr} comes from the fact that the
formulas for the homotopies given there are no longer applicable in
this case.  Let $\mu:S^1\times \RR\to S^1$ be given by
$\mu(\tau,t)=e^{2\pi i t} \tau$.  We also use the following convention
:
$d\tau\in\Omega^1S^1$ is the standard volume form on $S^1$ and for an
$S^1$-valued function $g$ we let $dg=g^{\ast}(d\tau)$.\\
(i) Let $F':M\times [0,1]\to S^1$ be the homotopy of $f_0$ and $f_1$.
By a simple partition of unity argument we may assume that $F'$ is flat
at the ends of the cobordism in the sense that there are collared
neighbourhoods $U_i\approx M\times [a_i,a_i+\delta]$ of $M\times\{i\}$,
$i\in\{0,1\}$, $a_0=0$, $a_1=1-\delta$, $\delta$ small, such that for
$(x,t)\in U_i$ we have $F'(x,t)=f_i(x)$.  Because $M$ is compact there
exists $m$ such that $|(\partial F'/\partial t)(x,t)|\leq m$ for all
$(x,t)\in M\times [0,1]$ (where $\partial F'/\partial t=\langle
dF',\partial/\partial t\rangle$).  Let $k:(-1/2,3/2)\to \RR$ be a Morse
function with exactly two critical points, a minimum at $0$ and a
maximum at $1$ that has value $1$ and such that for $t\in [\delta,1-
\delta]$ we have
$k'(t) < -m$.  Now define $F:M\times [0,1]\to S^1$ by $F=\mu (F'\times
k)$.  As $F'$ is flat close to $\partial (M\times [0,1])$, to verify
that $F$ is a Morse cobordism we only need to notice that $F$ has no
critical points in $M\times [0,1]\backslash (U_0\cup U_1)$.  This
happens because in this set $\partial F/\partial t = \partial
F'/\partial t +k'<0$ (which is computed by using polar
coordinates $\tau \to e^{2\pi i\tau}$
on $S^{1}$).  It is obvious that there are metrics $\alpha'$
that extend $\alpha_i$ and such that $(F,\alpha')$ is Morse-Smale.\\
(ii) This is the analogue of the construction of linear Morse
cobordisms in Lemma \ref{constr3}.  We consider a Morse function
$k:(-1/2,3/2)\to \RR$ as above (without any derivative restrictions).
Let $\nu : [0,1]\to [0,1]$ be an increasing $C^{\infty}$-function with
$\nu([0,\delta])=0$, $\nu([1-\delta,1])=1$.  We let $F':M\times
[0,1]\to S^1$ be a smooth homotopy of $f_0$ and $f_1$ that is
defined by $F'(x,t)=e^{2\pi i\nu(t)}f_0(x)$.  Finally, we define $F$
as before $F(x,t)=\mu(F'\times k)$.  We take on $M\times [0,1]$ the
metric $\alpha'=\alpha_0+dt^2$.  It is easy to verify that
$(F,\beta)$ is a cobordism whose induced morphism is the identity.
Indeed, the lift of this pair to $\overline{F}:\overline{M}\times
[0,1]\to \RR$ is an obvious linear cobordism in the real valued sense
and this implies that $F^{M\times [0,1]}$ is the identity.
\end{proof}

\noindent
(i) Given these two lemmas item (i) of Theorem \ref{theo:Novikov}
follows exactly by the
same argument as item (ii) of Theorem \ref{theo:rig}.  More precisely,
assume that $(f_i,\alpha_i):M\to S^1$ are homotopic Morse-Smale
functions such that $f_1$ is $C^0$ close to $f_0$.  Then there exists a
homotopy $F':M\times [0,1]\to S^1$ of $f_0$ to $f_1$ which is
$C^0$-close to the constant homotopy $f_0$.  This means that
by taking $f_1$ sufficiently $C^0$ close to $f_0$ we may also assume
the constant $c$ in Lemma \ref{constr4} (i) to be as small as desired.
We now construct Morse cobordisms relating $e^{2\pi i c'}f_0$, $e^{2\pi
i c}f_1$ and $f_0$ for small constants $c$ and $c'$ that are similar to
those in Lemma \ref{two-parameter} and the argument in Lemma
\ref{setting}
applied to the lifts of $\overline{f}_0$, $\overline{f}_1$ to
$\overline{M}$ implies that the resulting chain maps
$$\begin{array}{l}
w~:~C^{Nov}(M,f_0,\alpha_0)\to C^{Nov}(M,f_1,\alpha_1)~,\\[1ex]
w'~:~C^{Nov}(M,f_1,\alpha_1)\to C^{Nov}(M,f_0,\alpha_0)
\end{array}$$
have the property
that $w'\circ w$ satisfies the assumptions of Lemma \ref{lem:algNov}
and is therefore an isomorphism.  If $f_0=f_1$ the fact that $w'$ is
surjective implies that it is an isomorphism (because any module
epimorphism between finitely generated free modules of the same rank is
an isomorphism).

\noindent
(ii) Again, the proof here consists of adapting the proof of Theorem
\ref{theo:glueadapt} (i) to the $S^1$-valued case.  We shall use the
notations introduced in Definition \ref{splitting2}.  In particular $f$
and $h$ only differ in the interior of a small tubular neighbourhood
$U=N\times [-\epsilon, \epsilon]$ and on this set, if we identify the
interval $[-\epsilon,\epsilon]$ with a small arc we have $f(x,t)=t$.
By taking $\tau$ small we can make the function $h|_{U}$ as $C^0$ close
to the function $f|_{U}$ as desired.  Therefore, we may find a small
constant $l$ such that $h(x,t)<t+l=f(x,t)+l$.  We then construct by the
methods in Lemma \ref{constr3} a Morse function $F^1:U\times [0,1]\to
[-\epsilon,\epsilon +c]$ which is of the form $(1-u(t))h+u(t)(f+l)$
with $u(1)=1$, $u(0)=0$, $u'(0)=u'(1)=0$, $u''(0)>0>u''(1)$ and such
that $(\partial F^1/\partial t)(y,t) <0$ if $(y,t)\in U\times (0,1)$.
Because, $h$ and $f$ coincide outside $U$ we notice that on $\partial
U$ we have $F^1(y,t)=f(y)+lu(t)$.  Recalling that actually the image of
$F^1$ is in $S^1$ we rewrite $F^1|_{\partial U\times [0,1]}$ as
$F^1(y,t)= e^{2\pi i l u(t)}f$.  This shows that $F^1$ extends to a
function $H:M\times [0,1]\to S^1$ which is equal to $e^{2\pi i
lu(t)}f(y)$ if $(y,t)\in (M\backslash U)\times [0,1]$.  With an
appropriate metric $\alpha'$ the function $H$ is a Morse cobordism
between $(h,\alpha)$ and $(e^{2\pi i l}f,\alpha)$.  By inspecting its
lift
$\overline{H}$ we see that the induced morphism $T=H^{M\times
[0,1]}:C^{Nov}(M,f,\alpha)\to
C^{Nov}(M,h,\alpha)$
has the following property: for each generator $x\in {\rm
Crit}_{\ast}(f)$ we have $T(x)=x+o(x)+zp$ with $o(x)\in \ZZ[\pi]\{{\rm
Crit}_{\ast}(h)\backslash {\rm Crit}_{\ast}(f)\}$ and $p\in
\ZZ[\pi]_{\zeta}[[z]]\{{\rm Crit}_{\ast}(h)\}$.  We consider the
composition $J'=J_h\circ T$.  The description of $J_h$ now shows that
$J'(x)=x+zq$ with $q\in\ZZ[\pi]_{\zeta}[[z]]\{{\rm Crit}_{\ast}(f)\}$.
Lemma \ref{lem:algNov} implies that $J'$ is an isomorphism.

\noindent
(iii) We intend to prove this item by showing that for a fixed $n$ and
some special choice of function $g_n$ leads by the construction above
to a splitting $h_n$, a cobordism $(H_n,\alpha_n)$ and a corresponding
morphism
$T_n$ that has the form $T_n(x)=x+o_n(x)+z^np$ such that with the
notations in Definition \ref{splitting2} we have
$$p\in \ZZ[\pi]_{\zeta}[[z]]\{{\rm Crit}_{\ast}(h)\}$$
and
$$o_n\in\ZZ[\pi]_{\zeta}[[z]]\{{\rm Crit}_{\ast-1}(g_n)
\times \{-\epsilon/2\}\}~.$$ Notice that if such a $T_n$ is
constructed, then as $J_{h_n}:C^{Nov}(M,h_n)\to \widehat{F(h_n)}$
vanishes on $\ZZ[\pi]_{\zeta}[[z]]\{{\rm Crit}_{\ast-1}(g_n)\times
\{-\epsilon/2\}\}$ the composition $J_{h_n}\circ T_n$ is an isomorphism
of the form required.  Therefore, we have reduced the proof to the
construction of $h_n$, $H_n$ such that $T_n$ is as above.  This
construction may be accomplished in a way similar to the proof of
Theorem \ref{theo:glueadapt} (iii).  We recall that the metric $\alpha$
is the product metric $\beta +dt^2$ inside a neighbourhood $N\times
[-\delta,\delta]$.\\
\indent As in Theorem \ref{theo:glueadapt} (iii), the key point
is to work with an adapted function $g$ (see Definition \ref{adapted}).
In our setting, consider a Morse-Smale function $(g_n,\beta):N\to \RR$
which is adapted to splitting the Morse-Smale function
$(f_N(i),\alpha_N(i)):M_N(i)\to \RR$ along $z^{i-1}N$ for all
$0<i\leq n$ (the notations are as in Lemma \ref{lem:Nov}).  Such a
function exists by the same argument as in the real valued case (there
are more conditions to be verified but they are of the same type and
finite in number).  We now construct the $(\epsilon,\tau)$-splitting
$h_n$
such that inside the set $N\times [-\epsilon,\epsilon]$ the function
$h_n$ has the particular form of the functions $f_{\epsilon}$ in the
proof of Proposition \ref{theta1}.  There is one minor additional point
that needs to be looked at: we need to be sure that $(h_n,\alpha)$
satisfies the Morse-Smale transversality condition.  The needed
transversality is certainly satisfied with respect to the unstable and
stable manifolds of the critical points in ${\rm Crit}_{\ast}(h)\cup
{\rm Crit}_{\ast}(f)$ by the transversality part in the ``adapted"
condition.  The rest of the transversalities needed can be insured by
possibly slightly perturbing the {\em stable} manifolds of $g_n$.  This
can preserve the transversalities already achieved and does not affect
the part (i) of the adapted condition (see Definition \ref{adapted}).\\
\indent
We now construct the Morse cobordism $(H_n,\alpha_n):M\times [0,1]\to
S^1$ in the
same general way as for item (ii), but now we need to control the
behaviour of the flow lines of $\overline{H}_n$ over the set
$K=M_N(n)\backslash (z^n+z^{-1})(N\times (-\epsilon,+\epsilon))$.
More precisely, we shall construct a cobordism $(H_n,\alpha_n)$ which
satisfies
\begin{equation}\label{eq:van_fl} n^{\overline{H}_n}(x,y)=0
\end{equation}  for  all $x\in {\rm Crit}_k(f)\times \{1\}\subset
\overline{M}\times \{1\}$ and
$$y\in z^i({\rm Crit}_k(f)\cup ({\rm Crit}_k(g_n)\times
\{\epsilon/2\}))$$
for all $0<i\leq n-1$.  This immediately implies our claim for $T_n$
(as in the proof of Theorem \ref{theo:glueadapt}).  The argument
providing this $(H_n,\alpha_n)$ follows.  It is quite similar to that
in the proof of Theorem \ref{theo:glueadapt} (iii), but with the
additional difficulty that all the constructions initially accomplished
for real valued functions (on $\overline{M}$) will need to be adjusted
such that they are periodic in order to descend to $M$.\\
\indent
We first make the properties of $h_n$ more explicit.  We assume that
the neighbourhood $N\times [-2\epsilon,2\epsilon]$ is like in the
definition of splittings (after reparametrization $f$ has the form
$(x,t)\mapsto t$ and the metric is the product one here).  On $N\times
[-2\epsilon,2\epsilon]$ we consider a Morse-Smale function $g'$ such
that $g'(x,t)=t$ for $(x,t)\in N\times ([-2\epsilon,-\epsilon]\cup
[\epsilon,2\epsilon])$; $g'$ is a linear cobordism between
$\epsilon(g_n+k)$ and $\epsilon(g_n+k')$ on $N\times
[-\epsilon/2,\epsilon /2]$ and is a linear combination of the form
$u(t)(-\epsilon)+(1-u(t))\epsilon(g_n+k')$ inside $N\times
[-\epsilon,-\epsilon/2]$ with $u:[-\epsilon,-\epsilon/2]\to \RR$ an
appropriate function with $u([-2\epsilon -\epsilon)=1$ and
$u(-\epsilon/2)=0$ as in Lemma \ref{constr1}; a similar linear
combination is valid inside $N\times [\epsilon/2,\epsilon]$.  The
function $h_n$ equals $g'$ on $N\times [-\epsilon,\epsilon]$ and is
equal to $f$ on $M\backslash N\times [-\epsilon,\epsilon]$.\\
\indent
Similarly, we may define a standard cobordism $(H',\alpha')$ between
$g'$ and the function $(x,t)\mapsto t+c$ which is defined on $(N\times
[-2\epsilon,2\epsilon])\times [0,1]$ and where $c$ is a very small
constant depending on $k,k',\epsilon$.  This cobordism $H'$ differs
from a linear cobordism only in the interior of $N\times
[-\epsilon,\epsilon]\times [0,\epsilon]$.  The metric $\alpha'$ is the
product metric $\alpha +d\tau^{2}$.  We fix $\delta>0$ such that if
$(x,t,t')\in N\times [-2\epsilon,2\epsilon]\times [1-\delta,1]$, then
the flow line of $-\nabla^{\alpha'}H'$ containing $(x,t,t')$ does not
intersect $N\times [-2\epsilon,2\epsilon]\times [0,2\epsilon]$.\\
\indent
Now let $0<\epsilon_i<\epsilon<1/3 $ and consider $\phi_i:N\times
[-2\epsilon,2\epsilon ]\times [0,2\epsilon]\to N\times
[-2\epsilon,2\epsilon]\times [0,2\epsilon]$ a diffeomorphism given by
$(x,(t,t'))\mapsto (x,l_i(t,t'))$ where
$l_i:[-2\epsilon,2\epsilon]\times [-2\epsilon,2\epsilon]\to
[-2\epsilon,2\epsilon]\times [-2\epsilon,2\epsilon]$ is a
diffeomorphism such that $l_i$ restricts to the identity on the
complement of the square of side $[-(3/2)\epsilon, (3/2)\epsilon]$ and
the image of the square of side $[-\epsilon_i,\epsilon_i]$ is the
square of side $[-\epsilon,\epsilon]$.  We also require that $l_i$
restricts to diffeomorphisms when $t=0$ and when $t'=0$.  Let
$(H_i,h_i)= (\epsilon_i/\epsilon )(H', h')\circ \phi_i$.  The metric
associated to the Morse cobordism $H_i$ is the metric
$\alpha^i=\phi_i^{\ast}\alpha'$.  It is useful to note that $H_i$
coincides with a linear cobordism outside of $N\times
[-\epsilon_i,\epsilon_i]\times [0,\epsilon_i]$ and that
$\nabla^{\alpha^i}H_i(x,t,t')$ and $\nabla^{\alpha'}H_i(x,t,t')$
have the same projection on $N$.  Clearly, $\alpha^i$ and $\alpha'$
differ only outside $N\times [-(3/2)\epsilon,(3/2)\epsilon]\times
[0,(3/2)\epsilon]$.  For $0\leq i\leq n$ we fix the following notation
$K_i=z^i(N\times [-(3/2)\epsilon,(3/2)\epsilon])$ and $K'=K\backslash
\cup K_i$.

We intend to construct a cobordism $(H'',\alpha''):K\times [0,1]\to
\RR$ between a Morse-Smale function $h''$ and $f_N(n)$ such that there
exists a sequence of $\epsilon_i$'s with the property that $H''$ is a
linear cobordism on $K'\times [0,1]\times [0,1]$, on
$(H'',\alpha'')|_{z^i(N\times [-\epsilon,\epsilon] \times [0,1])}$ it
agrees with $(z^iH_i,\alpha^i)$ and it coincides with $z^iH'$ on
each set $z^i(N\times ([-2\epsilon,-(3/2)\epsilon]\cup
[(3/2)\epsilon,2\epsilon])\times [0,1])$ .  Moreover, $(H'',\alpha'')$
should verify the obvious analogue of (\ref{eq:van_fl}).\\
\indent
To construct $H''$ we proceed by induction.  We start with the trivial
linear cobordism with a product metric.  Assume now that we have
already constructed $(\epsilon_{n-1},\ldots \epsilon_i)$ and the
function $h'_i$ and $(H'_i,\alpha'_i)$ such that $H'_i|_{K_{j}\times
[0,1]}=z^jH_j$ for all $j<i$ and $H'_i$ is a cobordism between
$\overline{f}$ and $h'_i$ which is linear on $(K\backslash
\cup_{j<i}K_j)\times [0,1]$ and is such that the equation
(\ref{eq:van_fl}) is satisfied by $(H'_i,\alpha'_i)$.  We now want to
construct $\epsilon_{i-1},H'_{i-1},\alpha'_{i-1}$.\\
\indent
We first pick some choice for $R_{\epsilon'}=H'_{i-1}$ defined with
respect to some $\epsilon_{i-1}=\epsilon'$ and which satisfies the
properties required from $H'_{i-1}$ except possibly (\ref{eq:van_fl})
and such that $R_{\epsilon'}|_{K_{i-1}\times [0,1]}=z^{i-1}H_{i-1}$ and
$R_{\epsilon'}=H'_i$ in the complement of $K_{i-1}\times [0,1]$.  As in
the Morse case we might need to modify the associated metric on
$K\times [1/3,1-\delta]$ (in the real case it was enough to take
$\delta=3/2$ but here it is useful to use the constant fixed above)
such that $R_{\epsilon'}$ is Morse-Smale (see above Lemma
\ref{deform}).  In fact, it is easy to see that we may assume that the
metric $\alpha_{\epsilon'}$ associated to $R_{\epsilon'}$ satisfies
property (\ref{eq:van_fl}) with respect to $H'_i$ (because it can be
chosen to be arbitrarily close to $\alpha'_i$) in the exterior of
$K_{i-1}\times [0,2\epsilon])$), it coincides with the standard product
metric outside $K\times [1/3,1-\delta]\cup \cup_iK_i$, it coincides
with $\phi^{\ast}_i(\alpha')$ on $z^{j}N\times
[-2\epsilon,2\epsilon]\times [0,1]$, $j\leq i-1$, and is arbitrarily
close to the product metric $\alpha'$.  In fact, the set of choices for
the metrics $\alpha_{\epsilon'}$ is dense and open in some neighborhood
of $\alpha'$ viewed as a metric on $K'\times [0,1]$.  As the difference
between $R_{\epsilon'}$ and $H'_i$ appears around $z^{i-1}N$, as in the
proof of Lemma \ref{deform}, we assume that there exists a pair of
critical points $x,y$ of $R_{\epsilon'}$ that do not satisfy
(\ref{eq:van_fl}) even if we make $\epsilon'\to 0$.  By using the
induction hypothesis and proceeding as in the proof of that lemma (see
also Remark \ref{rem:class} b.) we are lead to the existence of a
sequence of points $x_1, x_2,\dots,x_t$ with the following properties:
$x_1$ is a critical point of $f_N(n)$ such that $x=x_1\times\{1\}$,
$$x_2,x_3\in \bigcup\limits_{{\rm ind}_{g_n}(q)<
{\rm ind}_{f}(x_1)} W^{u}_{g_n}(q)\subset z^{i-1}N~,$$
$x_t=y$, $x_3$ is related to $x_{4}$ by a flow line of $h'_i$, $x_{1}$
is related to $x_{2}$ by a flow line of $\overline{f}$, $x_{2}$ is
related to $x_{3}$ by a flow line of $g_n$ and all the points
$x_{4},\ldots, x_t$ are critical points of $h'_i$ and appear as
vertices in a sequence of broken flow lines of $h'_i$ (in the proof of
Lemma \ref{deform} this sequence only had 4 elements).  Moreover, and
this is a key point, all the points $x_{l}$ for $l>3$ belong to
$W^{u}_{H'_i}(x)$.  By transversality this means that ${\rm
ind}_{h'_i}(x_{l})\leq {\rm ind}(x_1)$.  By the transversality
properties of $g_n$ we have that
$${\rm ind}_{h'_i}(x_{4})<{\rm ind}_{f}(x_1)~=~{\rm ind}_{H'_i}(x)-
1~=~k-1~.$$
Notice that this implies that ${\rm ind}_{h'_i}(x_t)<{\rm ind}(x_1)$
which contradicts our assumption that ${\rm ind}(y)=k-1$.  Indeed, if
$t=4$ this has already been shown above and if $t>4$ then $y$ is in the
unstable manifold of the point $x_{t-1}$ whose index is at most $k-1$.
Again by transversality the index of $y$ is at most $k-2$.  Therefore,
the existence of the $\epsilon_i$ has been established and hence the
construction of $(H'',\alpha''),h''$ follows by recurrence.\\
\indent
In order to complete the proof we need to adjust $H'',h'',\alpha''$ in
such a way that the validity of equation (\ref{eq:van_fl}) is
preserved, but the new functions and metric are periodic and hence
induce corresponding notions on $M$.  Let $\psi:K\times [0,1]\to
K\times [0,1]$ be a diffeomorphism which leaves invariant $(K-\cup
z^i(N\times [-2\epsilon,2\epsilon]))\times [0,2\epsilon]$ and is the
inverse of $\phi_i$ on each set $z^i(N\times
[-2\epsilon,2\epsilon])\times [0,2\epsilon])$.  Consider the
composition $Q=H''\circ \psi$.  This function still satisfies
(\ref{eq:van_fl}) but, of course, with respect to the metric
$\alpha^{\ast}=\psi^{\ast}(\alpha'')$.  Notice that
$\alpha^{\ast}=\alpha'$ in $\cup_i( K_i\times [0,1]) \cup (K\times
[1-\delta ,1]\cup [0,1/3])$ and so it is immediate to see that we can
choose $\alpha''$ as described above on $K'\times [0,1]$ but also such
that $\alpha^{\ast}$ is periodic.  On each of the sets $z^i(N\times
[-2\epsilon,2\epsilon]\times [0,1])$ the function $Q$ satisfies
$Q(x,t,t')=u_i(t)z^iH'(x,t,t')$ where $u_i(t)$ equals $1$ on
$[-2\epsilon,-(3/2)\epsilon]\cup [(3/2)\epsilon,2\epsilon]$ and equals
$\epsilon_i/\epsilon$ on $[-\epsilon,\epsilon]$ (these functions
$u_i$ appear when the cobordisms $H_i$ are pasted with the linear
cobordisms on $K'\times [0,1]$).  This means that there exists a
diffeomorphism $\nu:[-n,1]\to [-n,1]$ which satisfies $\nu(t)=t$ for
$t\not\in \cup [k-2\epsilon,k+2\epsilon]$ and with the property that
for $(x,t,t')\in (N\times[k-2\epsilon,k+2\epsilon] \times [0,1])$, we
have $\nu\circ Q(x,t,t')=H'(x,t,t')+k$ (we have identified everywhere
$z^if(x)$ and $f(x)+i$).  We now define $H'''=\nu\circ Q$.  By
factorization this cobordism provides the desired cobordism
$H_n:M\times [0,1]\to S^1$.  Indeed, by construction, $\overline{H}_n$
satisfies the required properties relative to equation
(\ref{eq:van_fl}) with respect to the metric $\alpha^{\ast}$.  The pair
$(H_n,\alpha^{\ast})$ might not be a Morse-Smale pair even if the
Morse-Smale condition is satisfied up to order $n$.  Therefore, we
obtain the metric $\alpha_n$ by a new (and last) perturbation of
$\alpha^{\ast}$ away from $M\times \{0,1\}$ such that $\alpha_n$ is
sufficiently close to $\alpha^{\ast}$ for condition (\ref{eq:van_fl})
to continue to be satisfied (this is possible because this condition
concerns only the behavior of flow lines of ``length" $n$).  On
$M\times \{0,1\}$ the metric $\alpha_n$ coincides with $\alpha$.
Moreover, $H_n|_{M\times \{1\}}$ is equal to $f$ up to a constant and
$H_n|_{M\times\{0\}}$ is a splitting of $f$.  By construction the
morphism induced by $H_n$ satisfies equation (\ref{eq:van_fl}), and
this
concludes the proof.
\end{proof}

\begin{rem}\label{rem:fin1} {\rm
(a) For any fixed Morse-Smale function
$(f,\alpha):M\to S^1$, a splitting $(h,\alpha)$ of $(f,\alpha)$ along
$g:N=f^{-1}(0) \to \RR$ can be viewed as a geometric approximation of
finite nature for $(f,\alpha)$.  Indeed, $(f,\alpha)$ and $(h,\alpha)$
are compared by a simple Morse cobordism and the flow lines of $h$ are
never longer than one fundamental domain.  Theorem \ref{theo:Novikov}
(ii) indicates that the isomorphism type of the Novikov complex of
$(f,\alpha)$ can be recovered from any such finite approximation $h$
that is sufficiently $C^0$-close to $f$.  It is important to note that
$g$ can be here any Morse-Smale function $g:N\to \RR$.  Theorem
\ref{theo:Novikov} (iii) shows that by using $\ell$-adapted $g$'s in
the sense of being adapted to splitting $f_N(\ell)$ along $z^{\ell-1}N$
for higher and higher $\ell$'s we reconstruct the Novikov complex
itself.  Producing $\ell$-adapted $g$'s is a finite process that only
depends on $f_N(\ell)$ and its difficulty increases with $\ell$.  As in
any approximation process more effort is needed to produce better
approximations.\\
(b) A natural question is whether there is some $g:N\to \RR$ and an
associated splitting $h(g)$ such that there exists an isomorphism
$$\Theta^{\infty}:C^{Nov}(M,f,\alpha)\to \widehat{F(h(g))}$$ which is
basis preserving.  The major difficulty in proving such a statement
along the method above comes from the fact that one would need $g$ to
be adapted to splitting $f_N(\ell):M_N(\ell)\to \RR$ along $z^{\ell-
1}N$ for
any $\ell$.  The existence of such a $g$ is not at all clear.  It is
also
possible to use a slightly less restrictive condition at the point (ii)
in Definition \ref{adapted} but as it turns out this leads to some
technical difficulties related to transversality issues (in the
$S^1$-valued case).\\
(c)  The work of Pajitnov \cite{P} shows that for a
certain class $\mathcal{L}$ of pairs $(f,\alpha)$ which is
$C^0$-dense and open in the set of Morse-Smale pairs, there exists
an algebraic cobordism $(F,D,\theta,\theta',\psi)$ such that
the differentials of the Novikov complex of $(f,\alpha)$ are described
by the formula
$$\begin{array}{l}
d_{C^{Nov}(M,f,\alpha)}~=~
d_F+\sum\limits^{\infty}_{j=0}z^{j+1}\theta'\psi^j\theta~:\\[1ex]
C^{Nov}(M,f,\alpha)_i~=~(F_i)_{\zeta}((z)) \to
C^{Nov}(M,f,\alpha)_{i-1}~=~(F_{i-1})_{\zeta}((z))~.
\end{array}$$
with $\psi:D_i \to zD_i$ constructed by means of a homological
version of the partially defined flow-return.
Although $(F,D,\theta,\theta')$ was not constructed from a splitting
it certainly suggests that, at least for this class of functions,
the answer to the question
discussed in (b) is positive.  This seems even more likely
because the property that distinguishes $(f,\alpha)\in\mathcal{L}$ from
the rest of Morse-Smale functions appears to be quite close to
demanding the existence of an ``infinitely" adapted $g$.
It should be noted that it is a rather difficult task
to decide whether a particular Morse-Smale function
$(f,\alpha):M \to S^1$ belongs to $\mathcal{L}$ or not.  Therefore,
Theorem \ref{theo:Novikov} is a result essentially complementary to the
results in \cite{P} because \ref{theo:Novikov} applies to {\em all}
Morse-Smale
functions $(f,\alpha)$ the only choices involved having to do with the
auxiliary function $g$.\\
(d) As in Farber and Ranicki \cite{FR} define
the Cohn localization
$\Sigma^{-1}\ZZ[\pi]_{\zeta}[z,z^{-1}]$ of
$\ZZ[\pi]_{\zeta}[z,z^{-1}]$ inverting the set
$\Sigma$ of square matrices over
$\ZZ[\pi]_{\zeta}[z]$ sent to invertible
matrices over $\ZZ[\pi]$ by the augmentation $z
\mapsto 0$. Let
$(F,D,\theta,\theta',\psi)$
be the algebraic cobordism over $\ZZ[\pi]$
determined by a splitting $(h,\alpha)$ of
$(f,\alpha)$ as in Theorem \ref{theo:Novikov},
with
$$D~=~C(N,g,\beta)~,~F~=~C(M_N,f_N,\alpha_N)~.$$
The rational Novikov complex
$C^{FR}(M,f,\alpha)$ of \cite{FR} is the based
f.g. free
$\Sigma^{-1}\ZZ[\pi]_{\zeta}[z,z^{-1}]$-module
chain complex defined algebraically by
$$C^{FR}(M,f,\alpha)~=~
(
\Sigma^{-1}(F_i)_\zeta[z,z^{-1}],d_F+z\theta'(1-z\psi)^{-1}\theta)~.$$
The rational analogue of Theorem
\ref{theo:Novikov} shows that the
basis-preserving isomorphism type of
$C^{FR}(M,f,\alpha)$ is independent of the
choice of generic metric $\alpha$, and gives an isomorphism
$$\ZZ[\pi]_\zeta((z))\otimes_{\Sigma^{-1}\ZZ[\pi]_\zeta[z,z^{-1}]}
C^{FR}(M,f,\alpha)~\cong~C^{Nov}(M,f,\alpha)~.$$
See \cite{R} for further details of this
isomorphism, and of the construction of
$C^{FR}(M,f,\alpha)$. \hfill\qed}
\end{rem}

\providecommand{\bysame}{\leavevmode\hbox
to3em{\hrulefill}\thinspace}

\end{document}